\documentclass[smallextended]{svjour3}     

\smartqed  
\usepackage{graphicx}

\usepackage{savesym}
\usepackage{amsmath}
\usepackage{amssymb}
\usepackage{txfonts}      

\RequirePackage{natbib}
\usepackage{url}
\newcommand{\href}[2]{{#2}}

\RequirePackage{ifthen}

\newlength{\ei}\ei=0.0138888889em
\newlength{\eN}
\newlength{\SyW}  %
\newlength{\msu}  \msu=\mathsurround 
    
 \newcommand{\Ts}{\textstyle}
  \newcommand{\SSs}{\scriptscriptstyle}
\newcommand{\ssr}{\rm\scriptscriptstyle}
\newcommand{\req}{\relax}
\newcommand{\rfia}[1]{\makebox[\parindent][l]{%
                     \makebox[0em][r]{\rm(}\sf#1\rm)}}
\newcounter{ABCcB}
\newcommand{\theABCcC}{\alph{ABCcB}}

\newcommand{\Cov}{\mathop{\rm {{}Cov{}}}\nolimits} 
\newcommand{\Ew}{\mathop{\rm {{}E{}}}\nolimits} 
\newcommand{\Var}{\mathop{\rm Var}\nolimits}    
\newcommand{\ve}{\varepsilon}

\newcommand{\Lo}{\mathop{\rm {{}o{}}}\nolimits}
\newcommand{\LO}{\mathop{\rm {{}O{}}}\nolimits}

\newcommand{\tr}{\mathop{\rm{} tr{}}}
\newcommand{\B}{\mathbb B}

\newcommand{\R}{\mathbb R}

\newcommand{\N}{\mathbb N}
\newcommand{\Jc}{\mathop{\bf\rm{{}I{}}}\nolimits}

\newcommand{\iid}{\mathrel{\stackrel{\ssr i.i.d.}{\sim}}}

\newcommand{\Tfrac}[2]{\textstyle\frac{#1}{#2}}

\newcommand{\Tsum}{\mathop{\Ts\sum}\nolimits}

\newtheorem{Thm}{Theorem}
\newtheorem{Prop}[Thm]{Proposition}
\newtheorem{Lem}[Thm]{Lemma}
\newtheorem{Rem}[Thm]{Remark}
\newtheorem{Cor}[Thm]{Corollary}
\newtheorem{Def}[Thm]{Definition}

\newtheorem{Bez}[Thm]{Notation}

\numberwithin{equation}{section}
\numberwithin{Thm}{section}

\newcounter{ABCc}
\renewcommand{\theABCc}{\alph{ABCc}}
\newenvironment{ABC}{\begin{list}{
  \rfia{\theABCc}}{\usecounter{ABCc} \topsep 0ex \partopsep 0ex \itemsep0ex
  \parsep=\parskip \leftmargin 0em \rightmargin 0em \itemindent=\parindent
  \listparindent=\parindent  \labelsep 0.2em \labelwidth 0.5em }}{\end{list}}

\newcommand{\asy}{\mbox{as\hskip-4\ei .\hskip-16\ei.}\hskip12\ei}

\ifx\blinded\undefined
\title{Consequences of Higher Order Asymptotics for the MSE of M-estimators on Neighborhoods}

\author{Peter Ruckdeschel}

\institute{P. Ruckdeschel \at
              Fraunhofer ITWM, Department of Financial Mathematics, \\
              Fraunhofer-Platz 1, 67663 Kaiserslautern, Germany\\
              and Dept.\ of Mathematics, University of Kaiserslautern,\\
              P.O.Box 3049, 67653 Kaiserslautern, Germany \\
              \email{peter.ruckdeschel@itwm.fraunhofer.de}\\           
}
\else
\title{Consequences of Higher Order Asymptotics for the MSE of M-estimators on Neighborhoods}
\author{}\institute{}
\fi
\date{Received: date / Accepted: date}
\begin{document}
\maketitle
\begin{abstract}
In \citet{Ruck:03c}, we derive an asymptotic
expansion  of the maximal mean squared error  (MSE)
of location M-estimators on suitably thinned out,
shrinking gross error neighborhoods.
In this paper, we compile several consequences of this result:
With the same techniques as used for the ${\rm MSE}$, we determine higher order expressions
for the risk based on over-/undershooting
probabilities as in \citet{Hu:68b} and \citet{Ri:80b}, respectively.
For the MSE problem, we tackle the problem of second order robust optimality:
In the symmetric case, we find the second order optimal scores again of Hampel
form, but to an $\LO(n^{-1/2})$-smaller clipping height $c$
than in first order asymptotics. This smaller $c$ improves  ${\rm MSE}$
only by $\LO(n^{-1})$.
For the case of unknown contamination radius we generalize the
minimax inefficiency introduced 
in \citet{R:K:R:08}
to our second order setup. 
Among all risk maximizing contaminations we determine a ``most innocent'' one.
This way we quantify the ``limits of detectability''
in \citet{Hu:97b}'s definition for the purposes of robustness.
\end{abstract}

\section{Motivation/introduction}\req
%
%
This paper takes up the central result of \citet{Ruck:03c}: a uniform
higher order expansion of the means squared error (MSE) of
location M-estimators on suitably
shrinking and thinned out neighborhoods $\tilde{\cal Q}_n(r;\ve_0)$, repeated as
Theorem~\ref{main} in this paper for easier reference.
It is of the following form
\begin{equation}\boldmath\label{typresult}
  \sup_{Q_n\in\tilde{\cal Q}_n(r;\ve_0)}n\, {\rm MSE}(S_n,Q_n)=r^2 \sup |\psi|^2+ \Ew \psi^2 + \Tfrac{r}{\sqrt{n}\,}\,A_1+
  \Tfrac{1}{n}\,A_2+\Lo(\Tfrac{1}{n})
\end{equation}
Here $S_n$ is an M-estimator to socres $\psi$, and $A_1$, $A_2$ are polynomials in the contamination radius $r$,
in $b=\sup |\psi|$,
and in the moment functions
$t\mapsto \Ew \psi_t^l$, $l=1,\ldots,4$ and their derivatives evaluated in $t=0$, and $\ve_0$ is the breakdown
point of $S_n$, i.e.\ $\ve_0=\sup|\psi| / (\sup \psi - \inf\psi)$.
We recognize 
that the speed of the convergence to the first order \asy value is one order faster in the ideal model. 

In this paper we present some ramifications of this theorem, but in particular
consider its consequences for higher order robust optimality.
\begin{Bez}
  \rm\small
For indices we start counting with $0$, so that terms of first-order asymptotics
  have an index $0$, second-order ones a $1$ and so on. Also we abbreviate first-order, second-order and third-order
  by f-o, s-o, t-o respectively, and we write f-o-o, s-o-o, and t-o-o for first, second, and third-order asymptotically
  optimal respectively.
\end{Bez}

In Theorem~\ref{finminmax}, we take up
the over- and undershooting probabilities
used as risk in \citet{Hu:68b} to determine a finite sample
minimax estimator of location. By means of a s-o expansion, we refine
the corresponding f-o translation by \citet{Ri:80b},
providing a closer link to finite sample optimality. 

The closed form expressions in \eqref{typresult}, in particular under
certain symmetry assumptions, allows us to tackle corresponding
(uniform) higher order optimality problems, so that we may check
whether \citet{Pfan:79a}'s catchword ``{\it First order efficiency implies second order efficiency\/}''
survives when passing to neighborhoods around the ideal model,
which---at least under symmetry---indeed (partially) holds.

In this setting, we see that Huber-type location M-estimators remain optimal
in second order sense, and we even may determine the s-o-o clipping height
 $c_1=c_1(r,n)$ which in fact is slightly lower ($\LO(n^{-1/2})$) than
the f-o-o one. So in fact we only retain the optimal class, not the actual
 optimal estimator from f-o optimality.

For situations where the radius is (partially) unknown, the concept of a
\textit{minimax radius} has been introduced and determined in \citet{R:K:R:08}:
A radius $r_0$ is determined such that the
(f-o) maximal inefficiency $\bar \rho(r')$ (as defined in \eqref{rhobardef})
is minimized in $r'=r_0$. 
We translate this to the s-o setup; the s-o results in the Gaussian location model show that neither $c_1(r_1,\cdot)$,
nor s-o minimax radius $r_1(\cdot)$ vary much in $n$ and that for all $n$,
s-o minimax inefficiency is always smaller than the corresponding f-o one. 

Asymptotics also helps to understand  which contaminations are (already) dangerous:
We determine the \textit{cniper contamination} as a  most innocent appearing
least favorable contamination, which is shown to form a saddlepoint together
with the f-o (s-o) optimal M-estimator.
It appears to be innocent, as it produces only ``outliers'' which are hardest to detect in some
sense specified in this section. 



\paragraph{Organization of the paper}
We start with the setup
of one dimensional location and recall the main theorem of \citet{Ruck:03c} in
section~\ref{setupsec}. This result is generalized to a over-/undershooting probability
 loss in section~\ref{conseqsec}.

Consequences of Theorem~\ref{main} as to higher order robust optimality are discussed in
section~\ref{conseq2o}.
%
%
As a (partial) explanation for the good, respectively excellent behavior of f-o-o, s-o-o and t-o-o procedures
as to numerically exact finite maximal MSE, we present an argument
based on a functional implicit function theorem
in section~\ref{ImplMinThmS}.
For decisions upon the procedure to take, only relative risk is relevant
which is discussed in some detail in subsection~\ref{relrisk}.
Section~\ref{ramifsec} then considers further supplementary results to Theorem~\ref{main}:
a s-o variant of the minimax radius and s-o cniper contaminations.
The proofs to the theorems and propositions of this paper are collected in  
 section~\ref{proofsec}.

\section{Setup} \label{setupsec}
\subsection{One-dimensional location}
We consider estimation of parameter $\theta$ in a one--dimensional location model, i.e.\
\begin{equation}\label{locmod}
X_i=\theta+v_i, \qquad v_i\stackrel{\SSs \rm i.i.d.}\sim F,\qquad   P_{\theta}={\cal L}(X_i)
\end{equation}
for some ideal distribution $F$ with finite Fisher-Information of location ${\cal I}(F)$,
i.e.\
\begin{equation}\label{FI}
\Lambda_f=-\dot f/f \in L_2(F),\qquad{\cal I}(F)=\Ew[\Lambda_f^2]<\infty
\end{equation}
We also assume that $\Lambda_f$ is increasing.
By translation equivariance, we may restrict ourselves to $\theta_0=0$ which will be suppressed in the notation.

The set of {\em influence curves\/} (IC's) $\Psi$ for the estimation of $\theta$ is defined
as \citet{Ri:94}
\begin{equation}
  \Psi:=\{\psi \in L_2(F) \,|\, \Ew[\psi]=0,\quad \Ew[\psi \Lambda_f]=1\},
\end{equation}
where both expectations are evaluated under $F$.
As class of estimators we consider {\em asymptotically linear estimators\/} (ALE's), i.e.\
estimators $S_n=S_n(X_1,\ldots,X_n)$ with the property
\begin{equation} \label{ALEd}
\sqrt{n\,} \,\,S_n=\Tfrac{1}{\sqrt{n}}\Tsum_{i=1}^n \psi(X_i)+\Lo_{F^n}(n^0)
\end{equation}
We consider maximal mean squared error (MSE) on shrinking neighborhoods
 of this ideal model,
defined as the set ${\cal Q}_n(r)$ of distributions
\begin{equation}\label{contadef}
{\cal L}^{\SSs \rm real}_{\theta}(X_1,\ldots,X_n)=Q_{n}=
{\Ts\bigotimes\limits_{i=1}^n}[(1-\Tfrac{r_n}{\sqrt{n}}) F+\Tfrac{r_n}{\sqrt{n}}\,P^{\SSs \rm di}_{n,i}]
\end{equation}
with $r_n=\min(r,\sqrt{n})$, $r>0$ the contamination radius and $P^{\SSs \rm di}_{n,i} \in {\cal M}_1(\B)$ arbitrary,
uncontrollable contaminating distributions.
As usual, we interpret $Q_n$ as the distribution
of the vector $(X_i)_{i\leq n}$ with components
\begin{equation}\label{Uidef}
  X_i:=(1-U_i)X_i^{\SSs \rm id}+U_i X_i^{\SSs \rm di},\qquad i=1,\ldots,n
\end{equation}
for $X_i^{\SSs \rm id}$, $U_i$, $X_i^{\SSs \rm di}$ stochastically independent, $X_i^{\SSs \rm id}\iid F$, $U_i\iid{\rm Bin}(1,r/\sqrt{n})$,
and $(X_i^{\SSs \rm di})\sim P^{\SSs \rm di}_{n}$ for some arbitrary $P^{\SSs \rm di}_{n} \in{\cal M}_1(\B^n)$.

Suppressing the dependency upon $\theta$ as usual, in \citet{Ri:94}, the first order expansion of maximal MSE
of an ALE is derived as
\begin{equation}\label{MSEpsi}
\tilde R(S_n,r)=r^2\sup |\psi|^2 + \Ew_{\rm \SSs id} |\psi|^2
\end{equation}
The (first-order) MSE-optimal IC $\eta_{b_0}$  
in a smooth $p$-dimensional parametric model with $L_2$-derivative $\Lambda$
by Theorem~5.5.7 (ibid.) has to be of Hampel form
\begin{equation} \label{allgemoptIC}
\eta_{b_0} = Y\min\{1,b_0/|Y|\},\qquad Y=A\Lambda-a
\end{equation}
for some $A\in\R^{p\times p}$, $a\in\R^p$ such that $\eta_{b_0}$ is an IC,
and  $b_0$ solving
  $\Ew(|Y|-b_0)_+=r^2b_0$.
In our location context, for Lagrange multipliers $z$ and $A$ such that $\eta_{b_0}=\eta_{c_0} \in \Psi$, we get that
\begin{eqnarray}\label{HK1}
&&\eta_{c_0}=A(\Lambda_f-z)\min\{1,
c_0/{|\Lambda_{f}-z|}\},\\
&&c_0 \;\;{\rm s.t.}\quad \Ew[(|\Lambda_f-z|-c_0)_+]=r^2c_0 \label{HK3}
\end{eqnarray}
\subsection{Higher Order Expansion} \label{HOE}
In \citet{Ruck:03c} we obtain corresponding higher order expansions of the maximal MSE
if we thin out the neighborhood system to the set $\tilde {\cal Q}_n(r;\ve_0)$
of conditional distributions
\begin{equation}\label{modifdefi}
Q_{n}={\cal L}\Big\{[(1-U_i)X_i^{\SSs \rm id}+U_i X_i^{\SSs \rm di}]_i\,\Big|\,\sum U_i \leq\,\ulcorner \ve_0 n\, \urcorner-1\,\Big\}
\end{equation}
where $\ve_0=1/({2+\delta_0})$ is the functional (\citet[(2.39),(2.40)]{Hu:81}) and the finite sample ($\ve$-contamination) breakdown
point (\citet[section 2.2]{Do:Hu:83}) of the corresponding M-estimator and $\delta_0$ is defined by
\begin{equation}\label{checkbdef}
\check b:= \inf\psi,\qquad\hat b =\sup\psi,\qquad  \bar b:=\Tfrac{1}{2}(\hat b -\check b),\qquad
\delta_0:=\Tfrac{|\check b+\hat b|}{\min((-\check b),\hat b)}\ge 0
\end{equation}
For the result we use the following assumptions and notation:
To scores function $\psi:\R\to\R$ let $\psi_t(x):=\psi(x-t)$ and define the following
functions
$L(t):=\Ew \psi_t$, $\psi_t^0:=\psi_t-L(t)$, $V(t)^2:=\Var \psi_t$,
$\rho(t):=\Ew(\psi_t^0)^3/V(t)^3$, $\kappa(t):=\Ew[(\psi_t^0)^4]/V(t)^4-3$.
Let $\check y_n$ and  $\hat y_n$ 
sequences in $\R$ such that for some $\gamma>1$,
$\psi(\check y_n)=\inf \psi +\Lo(\Tfrac{1}{n^\gamma})$, $\psi(\hat y_n)=\sup \psi +\Lo(\Tfrac{1}{n^\gamma})$.
For $H\in{\cal M}_1(\B^n)$ and an ordered set of indices $I=(1\leq i_1<\ldots<i_k\leq n)$ denote $H_{I}$ the marginal of
$H$ with respect to $I$. Consider three sequences $c_n$, $d_n$, and   $\kappa_n$ in $\R$, in $(0,\infty)$, and in $\{1,\ldots,n\}$, respectively.
We say that the sequence $(H^{(n)})\subset{\cal M}_1(\B^n)$ is {\em $\kappa_n$--concentrated left [right] of $c_n$ up to $\Lo(d_n)$}, if
for each sequence of ordered sets $I_n$ of cardinality $i_n\leq \kappa_n$
$
1-H^{(n)}_{I_n}\big((-\infty;c_n]^{i_n}\big)=\Lo(d_n)
$, $
\Big[\,
1-H^{(n)}_{I_n}\big((c_n,\infty)^{i_n}\big)=\Lo(d_n)\,
\Big]$.
For the theorem we make the following assumptions:
\begin{itemize}
  \item[(bmi)]  $\sup\|\psi\| =b <\infty$, $\psi$ monotone, $\psi \in \Psi$
  \item[(D)]   For some $\delta\in(0,1]$, $L$, $V$, $\rho$,   and  $\kappa$ as defined above allow the expansions\vspace{-.5ex}
  \end{itemize}
\begin{align}
 & L(t)\,=\,l_1t\!+\!\Tfrac{1}{2} l_2\,t^2\!+\!\Tfrac{1}{6}l_3\,t^3\!+\!\LO(t^{3+\delta}),&\hspace{-.5em}&
   V(t)\,=\,v_0(1\!+\!\tilde v_1\,t\!+\! \Tfrac{1}{2}\tilde v_2\,t^2)\!+\!\LO(t^{2+\delta})\label{VD}\\
&\rho(t)\,=\,\rho_0\!+\!\rho_1\, t\!+\!\LO(t^{1+\delta}),&\hspace{-.5em}&
  \kappa(t)\,=\,\kappa_0\!+\!\LO(t^{\delta})\label{KD}
\end{align}
\begin{itemize}
  \item[(Pd)]   There are some $T>0$ and $\eta>0$ such that
  \begin{equation}\label{abkling}
F(t)\geq 1-t^{-\eta},\quad \mbox{for } t>T,\qquad F(t)\leq (-t)^{-\eta}\quad\mbox{ for } t<-T
\end{equation}
  \item[(C)] Let $f_t$ be
  the characteristic function of $\psi_t(X^{\rm\SSs id})$; then
  \begin{equation}\label{(C)gl}
\lim_{t_0\to 0}\limsup_{s\to\infty} \sup_{|t|\leq t_0}|f_t(s)|<1
  \end{equation}
  \end{itemize}
%
With these preparations, we have the following theorem (\citet[Thm.~3.5]{Ruck:03c})
\begin{Thm}\label{main}
  In our one-dim.\ location model assume {(bmi)} to {(C)}  
\begin{ABC}
  \item the maximal MSE of the M-estimator $S_n$ to scores-function $\psi$
  expands to
  \begin{eqnarray}
    R_n(S_n,r,\ve_0)
    &=&{r}^{2}{b}^{2}+{v_0}^{2}+\Tfrac{r}{\sqrt{n}\,}\,A_1+\Tfrac{1}{n}\,A_2+\Lo(n^{-1})\label{mainres}
  \end{eqnarray}
with
\begin{small}
\begin{eqnarray}
A_1&=&  {v_0}^{2}\,\Big( \pm(4\,\tilde v_1+3\,l_2 \,)b+1 \Big)+{b}^{2} +
 [2\,{b}^{2}\pm l_2\,{b}^{3} \,]\,{r}^{2}\label{A1defr}\\
A_2&=&{{v_0}^{3}\,\Big((l_2+2\,\tilde v_1 \,)\rho_0+\Tfrac{2}{3}\,\rho_1\Big)+
 {v_0}^{4}\,(3\,\tilde v_2+{\Tfrac {15}{4}}\,{l_2^2}+l_3+9\,{\tilde v_1}^{2}+
 12\,\tilde v_1\,l_2 \,)}+\nonumber\\
 &&\quad+[\, {v_0}^{2}\,\Big( (3\,\tilde v_2+3\,{\tilde v_1}^{2}+\Tfrac{15}{2}\,{l_2^2}+2\,l_3+
 12\,\tilde v_1\,l_2 \,){b}^{2}+1\pm (8\,\tilde v_1+6\,l_2 \,)\,b \Big)+\nonumber\\
 &&\quad\pm 3\,l_2\,{b}^{3}+
 5 \,{b}^{2} \,]\,{r}^{2}+\Big( (\Tfrac{5}{4}\,{l_2^2}+\Tfrac{1}{3}\,l_3 \,){b}^{4}
 \pm 3\,l_2\,{b}^{3} +3\,{b}^{2}\Big)\,{r}^{4}\label{A2defr}
\end{eqnarray}
\end{small}
and we are in the $-\,[+]$-case depending on whether  \eqref{contbed1} or \eqref{contbed2} below applies.
\item
let $P^{\SSs \rm di}_{n}:=\bigotimes_{i=1}^n P^{\SSs \rm di}_{n,i}$ be contaminating measures for \eqref{contadef}.
Then $Q_n$ with $P^{\SSs \rm di}_{n}$ as  contaminating measures generates maximal risk in \eqref{mainres} if
for $k_1>1$ and $k_2>2\vee (\frac{3}{2}+\frac{3}{2\delta})$ with $\delta$ from (Vb)
and $K_1(n)=\ulcorner k_1 r\sqrt{n} \urcorner$
either
\begin{equation}
  \label{contbed1}
(P^{\SSs \rm di}_{n}) \mbox{ is  $K_1(n)$--concentrated left of $\check y_n- b\sqrt{k_2\log(n)/n}$ up to }\Lo(n^{-1})
\end{equation}
or
\begin{equation}
  \label{contbed2}
(P^{\SSs \rm di}_{n}) \mbox{ is  $K_1(n)$--concentrated right of $\hat y_n+ b\sqrt{k_2\log(n)/n}$ up to }\Lo(n^{-1})
\end{equation}
More precisely, if $\sup\psi <[>] -\inf \psi$, the maximal MSE is achieved by contaminations according to \eqref{contbed1}
[\eqref{contbed2}]. In case $\sup \psi=-\inf\psi$,  \eqref{contbed1} [\eqref{contbed2}] applies if
\begin{equation} \label{wokontlra}
\tilde v_1\,>\,[<]-\Tfrac{l_2}{4}\Big(\Tfrac{b^2}{v_0^2}(r^2+3)(1+\Tfrac{r}{\sqrt{n}}-\Tfrac{2r^2}{n})+3(1-\Tfrac{b^2}{v_0^2})\Big)
\end{equation}
If $\sup \psi=-\inf\psi$ and there is ``$=$'' in \eqref{wokontlra}, \eqref{contbed1} and \eqref{contbed2} generate the same
risk up to order $\Lo(n^{-1})$.
\end{ABC}
\end{Thm}
\paragraph{Special cases}
Let $Q_n^0$ be any distribution in $\tilde {\cal Q}_n$ attaining maximal risk
in Theorem~\ref{main}. Under symmetry or more specifically if
\begin{equation} \label{symdef}
l_2=v_1=\rho_0=0,
 \end{equation}
we obtain as maximal risk in \eqref{mainres}

\begin{small}
\begin{eqnarray}
n\,\Ew_{Q_n^0}[S_n^2\,]&=&
\left(  {r}^{2}{b}^{2}+{v_0}^{2}\right) \left(1+\Tfrac{r}{\sqrt{n}}+\Tfrac{r^2}{n}\right)+\Tfrac{r}{\sqrt{n}}\,
 \left(b^2(1+r^2)\right)+
\Tfrac{r^2}{n}\, \left(b^2(5+2r^2)\right)+{\frac {\frac{2}{3}{v_0}^{3}\,\rho_1+
 {v_0}^{4}\,(3\,\tilde v_2+l_3\,)}{n}}+\nonumber\\
 &&\quad+
 {\frac {\left( {v_0}^{2}\, (3\,\tilde v_2+2\,l_3\,)\,{b}^{2} \right)\,{r}^{2}
 +\frac{1}{3}\,l_3 {b}^{4}\,{r}^{4} }{n}}+\Lo(n^{-1}), \label{symEF}
\end{eqnarray}
\end{small}

while under $r=0$ (with or without \eqref{symdef}), we get

\begin{small}
\begin{equation}
n\,\Ew_{F^n}[S_n^2\,]=
 {v_0}^{2}+{\frac { {v_0}^{3}\,\left((l_2+2\,\tilde v_1 \,)\rho_0+\frac{2}{3}\,\rho_1\right)}{n}}+
 {\frac {{v_0}^{4}\,(3\,\tilde v_2+l_3+{\frac {15}{4}}\,{l_2^2}+
 12\,\tilde v_1\,l_2+9\,{\tilde v_1}^{2} \,)}{n}}+\Lo(n^{-1})
\end{equation}
\end{small}

respectively, again under \eqref{symdef},

\begin{equation} \label{symr=0}
n\,\Ew_{F^n}[S_n^2\,]=
 {v_0}^{2}+{\frac {\frac{2}{3}\,{v_0}^{3}\rho_1+
 {v_0}^{4}\,(3\,\tilde v_2+l_3\,)}{n}}+\Lo(n^{-1}).
\end{equation}

\section{Other loss functions}\label{conseqsec}\req
One easily shows that under similar condition as for Theorem~\ref{main}, we may replace the squared loss function in the MSE
 by other loss functions
$\ell$ growing atmost at a polynomial rate. In this respect, Theorem~\ref{main} easily extends to
uniform convergence of other risks on $\tilde {\cal Q}_n$, e.g.\ absolute error ($\ell(x)=|x|$), $L_k$-error ($\ell(x)=|x|^k$)
for $1<k<\infty$, and certain covering probabilities, $\ell(x)=\Jc_{(\alpha_1,\alpha_2)}(x)$ for some $\alpha_1<\alpha_2\in\R$.\\
As an illustration, we consider this last type of loss function, more specifically in the form in which it arises in the
finite minimax estimation theory as in \citet{Hu:68b} and in which it has been extended to an \asy setup by \citet{Ri:80b}:
The risk is defined as
\begin{equation}
R^\natural(S_n,r)=\sup_{Q_n \in {\cal Q}_n(r)} \max\{ Q_n(S_n>\theta+\frac{\alpha_2}{\sqrt{n}}),\,Q_n(S_n<\theta-\frac{\alpha_1}{\sqrt{n}})\}
\end{equation}
\citet{F:Y:Z:01}  have taken up a similar setup with conventional confidence intervals to cover bias and variance simultaneously.
We work in the setup of \citet{Ri:80b} here and confine ourselves to the higher order terms of order $n^{-1/2}$, but of course an
 extension to terms up to order $n^{-1}$ as in Theorem~\ref{main} is feasible. Due to translation equivariance, it is no restriction to
consider the case $\theta=0$ only. As in \citet{Ri:80b}, we work with
a possibly asymmetric partition of the interval of given length $2a/\sqrt{n}$ laid around the estimator: Using the partition
\begin{equation}\label{alph12}
  2a=\alpha_1+\alpha_2=\alpha_1(S_n)+\alpha_2(S_n),
\end{equation}
we minimize the risk according to \citet[formulas (2.8) and (2.11) in]{Ri:80b}, if
with $\check b$, $\hat b$, and $\bar b$ from \eqref{checkbdef}
and
\begin{equation}
  \alpha_1=a-\delta,\qquad \alpha_2=a+\delta,\qquad \delta=\Tfrac{r}{2}\,(\hat b+\check b)
\end{equation}
If we now account for terms of order $\Tfrac{1}{\sqrt{n}}$ we minimize the risk if we use the partition
\begin{equation}
  2a=\alpha_1'+\alpha_2'=\alpha_1'(S_n)+\alpha_2'(S_n),
\end{equation}
with
\begin{equation}
  \alpha_1'=a-\delta-\delta',\qquad \alpha_2=a+\delta+\delta',
\end{equation}
$\delta'=\delta'_n$ given in the theorem below.
To this end, let
\begin{equation} \label{s1def64}
s_1:=(-a+r \bar b)/{v_0}
\end{equation}
Then, with $\Phi$ and $\varphi$ c.d.f.\ and density of ${\cal N}(0,1)$ and using
the notation of Theorem~\ref{main}, we have
\begin{Thm}\label{finminmax}
 For the location model \eqref{locmod} of finite Fisher information \eqref{FI}, assume {(bmi)}, {(D')} and {(C')}.
 Then for sample size $n$,
 the minimal over-/undershooting probability of an M-estimator $S_n$ for scores-function $\psi$ in ${\cal Q}_n$
   obtains eventually in $n$ as
\begin{eqnarray}
R^\natural (S_n)&=&\!
  \sup_{Q_n \in {\cal Q}_n} \max\{Q_n(S_n\leq -\frac{\alpha_1'}{\sqrt{n}}),\,Q_n(S_n\geq \frac{\alpha_2'}{\sqrt{n}})\}=\nonumber\\
\!&=&\!  R_-(S_n,Q^0_{n;-})=R_+(S_n,Q^0_{n;+})
\end{eqnarray}
with $Q^0_{n;\,-}$ resp.\ $Q^0_{n;\,+}$ according to 
\eqref{contbed1} resp.\ \eqref{contbed2} and
\begin{eqnarray}
\lefteqn{R_-(S_n,Q^0_{n;-})=\Phi(s_1)+
\Tfrac{1}{\sqrt{n}\,v_0}\varphi(s_1)\times\nonumber}\\
&&\quad\times\Big[\Tfrac{ra}{2}+2{l_2 a \delta}-a s_1\tilde v_1 v_0- \Tfrac{r(\check b^2+\hat b^2)s_1}{4 v_0}+\Tfrac{r^2\bar b}{2}
\Big]+\Lo(\Tfrac{1}{\sqrt{n}})
\label{r1-}
\end{eqnarray}
and $\delta'=\delta'_n$ according to
\begin{equation}
\delta'=\Tfrac{1}{\sqrt{n}}\Big(-\Tfrac{r\delta}{2v_0}-\Tfrac{l_2}{2v_0}(a^2+\delta^2)-\tilde v_1 v_0s_1 \delta-
\Tfrac{\rho_0}{6} (s_1^2-1)+ \Tfrac{r\bar b \delta s_1}{v_0^2}+\Tfrac{r^2\delta}{2v_0}\Big)
\end{equation}
\end{Thm}
%
\begin{Rem}\rm\small
\begin{ABC}
  \item
  If $l_2=\tilde v _1=0$ and $\hat b=-\check b$, we obtain the same result as \eqref{r1-}, if we  use the expressions $b_n:={\rm Bias}_n$
  and $v_n^2={\rm Var}_n$ for bias and variance from \citet[Prop.6.4]{Ruck:03c}, plug them into the
  \asy risk, which gives  $\Phi((rb_n-a)/v_n)$, and then expand this up to $\Lo(n^{-1/2})$.
\item The numerical values obtainable by Theorem~\ref{finminmax}  should be compared to those of \citet[sections~11.3.3.3 and 11.4.1]{Koh:04di};
admittedly the approach of Theorem~\ref{finminmax} in this context gives rather poor (too liberal) approximations compared to those
in the cited reference (see the {\tt R}-file {\tt Thm31.R} available on the web-page to this article); this is
plausible though, as Kohl already starts with finitely optimal procedures whereas our approach improves upon
asymptotically optimal ones.
\end{ABC}
\end{Rem}

\section{Consequences: Higher Order Optimality and Relative Risk}\label{conseq2o}\req
In this section, we consider the class ${\cal S}_2$ of all M-estimators according
to  {(bmi)}, {(D')}, and {(C')} as well as {(Pd)}; 
correspondingly, we define ${\cal S}_3$ with {(D)}, {(C)} replacing {(D')}, {(C')};
we always assume that the class of M-estimators ${\cal H}$ of ICs of Hampel-type \eqref{HK1}
forms a subset of ${\cal S}_2$ [${\cal S}_3$]. In particular we assume $f$ to be log-concave.
\subsection{Second-order optimality}\label{soao}
Symmetry allows considerable simplifications;
for instance,  if $F$ is symmetric, i.e.\ $F(B)=F(-B)$ for all $B\in\B$, in \eqref{HK1}  always $z=0$.
But also, much deeper results are possible. Thus
for the rest of this subsection, we assume \eqref{symdef}. 
Then \eqref{symEF} gives the s-o-maximal MSE for any M-estimator
in ${\cal S}_2$; in particular
\begin{equation}\label{A1sym}
A_1=v_0^2+b^2(1+2 r^2)
\end{equation}
Condition~\eqref{symdef} clearly holds for skew symmetric $\psi$ and symmetric $F$.
For symmetric $F$, however, for any IC $\psi$, also $\tilde \psi:=-\psi(-\,\cdot\,)$
is an IC and hence so is the skew-symmetrized $\psi^{(\SSs s)}:=\frac{1}{2}(\psi+\tilde \psi)$, too.
But by convexity of the MSE, $\psi^{(\SSs s)}$ will be at least as good as $\psi$ as to MSE, hence
it is no restriction to only consider skew symmetric ICs, and we fall into the application range of
\citet[Thm.~3.1]{Ru:Ri:04}, i.e.,
\begin{Thm}\label{cvthm}
Assume that maximal \asy risk of an ALE on $\tilde {\cal Q}_n$ resp.\ $\tilde {\cal Q}'_n(\,\cdot\,,s_0)$
 is representable as $G(rb(\psi),v_0(\psi))$ for some convex real-valued function $G(w,s)$, strictly isotone in both arguments
 and totally differentiable, bounded away from the minimum for $w\to\infty$.
Then, on ${\cal Q}_n$, respectively on $\tilde {\cal Q}_n$,
the optimal IC\  of Hampel-type \eqref{HK1} for some clipping height $ b=Ac$ determined by
\begin{equation}\label{diffarg}
r\, v_0  \,  \partial_wG(rAc, v_0) = \partial_sG(r Ac, v_0)\,   A\Ew(| \Lambda- z|- c)_+
\end{equation}
\end{Thm}
In our case, this theorem specializes to
\begin{Cor}\label{cneqcor}
Assume a symmetric 
model~\eqref{locmod} with increasing $\Lambda_f$ and \eqref{FI}.
Under the assumptions of this section, the s-o-o M-estimator in class ${\cal S}_2$ has an IC of
 Hampel-type \eqref{HK1}  with $z=0$ and the s-o-o clipping height $c_1=c_1(n)$
is determined by
 \begin{equation}\label{saoc}
   r^2 c\,\Big(1+\frac{r^2+1}{r^2+r\sqrt{n}}\,\Big)=\Ew(|\Lambda|-c)_+
 \end{equation}
Always, $c_0>c_{1}(n)$. Suppose that 
$h(c):=\Ew(|\Lambda|-c)_+$ is differentiable in
$c_0$ with derivative $h'(c_0)$. Then,
\begin{equation} \label{c1def}
  c_1(n)=c_0\,\Big(1-\frac{1}{\sqrt{n}}\,\frac{r^3+r}{r^2-h'(c_0)}\,\Big)
  +\Lo(\frac{1}{\sqrt{n}})
\end{equation}
\end{Cor}
That is, (for $n$ large enough) {\bf the f-o-o clipping height $c_0$ always is too optimistic.}

Assume s-o risk of ICs of Hampel-type \eqref{HK1} is smooth enough in $c$ in its minimum~$c_1$ to allow a s-o Taylor expansion,
which is an assumption on the remainder $\Lo(n^{-1})$ present in \eqref{mainres}.
Then, around $c_1$, s-o risk behaves like a parabola. But,
as by \eqref{c1def}, $c_1-c_0=\LO(1/\!\sqrt{n}\,)$,  using $c_1$ instead of $c_0$
can only improve s-o risk by order $\LO(1\!/n)$. This even carries over to risks ``near'' s-o risk:
\subsection{Consequences for the exact MSE}\label{ImplMinThmS}
\begin{Prop}\label{ImplMinThm}
  Let $F, F_n, G_n\in {\cal C}_2(\R)$, $n\in\N$, such that for some $\beta\ge \beta'>0$
\begin{equation} \label{nahbei}
\begin{array}{ll}
{\rm (i)} &\sup_x |F_n-G_n| + |F'_n-G'_n| + |F''_n-G''_n|=\LO(n^{-\beta}),\\[2ex]
{\rm (ii)}&\sup_x |F_n-F| + |F'_n-F'| + |F''_n-F''|=\LO(n^{-\beta'})
\end{array}
\end{equation}
Assume that in $x_0\in\R$, $F(x_0)$ is minimal, and that $F''(x_0)=f_2>0$.
 Then
 \begin{ABC}
   \item there is some sequence $(x_n)\subset \R$ such that eventually in $n$, $F_n(x_n)$ is minimal and $\lim F''_n(x_n)=f_2$.
   \item  $|x_n-x_0|=\LO(n^{-\beta'})$.
\item there is some sequence $(y_n)\subset\R$ such that eventually in $n$,  $G_n(y_n)$ is minimal and $\lim_n G''_n(y_n)=f_2$.
\item $|y_n-x_n|=\LO(n^{-\beta})$.
\item $0\leq G_n(x_n)-G_n(y_n)=\LO(n^{-2\beta})$.
 \end{ABC}
\end{Prop}
%
%
The drawback of this proposition is that assumption~\eqref{nahbei} is difficult to check if we have
no explicit expression for $G_n$:
For given $r\ge 0$, let ${\rm asMSE}_{j=0,1,2}(c)$
be the f-o, s-o, and t-o maximal MSE
of an M-estimator in ${\cal H}$, 
and ${\rm exMSE}(c)$ the corresponding exact maximal MSE $R_n$; we would like to apply Proposition~\ref{ImplMinThm}
to $F={\rm asMSE}_0$, $F_n={\rm asMSE}_{j=1,2}$ 
and $G_n={\rm exMSE}$ to conclude on the performance of f-o-o, s-o-o, t-o-o procedures
as to ${\rm exMSE}$. 
As to \eqref{nahbei}, part (ii) is easy to see checking the expressions, giving $\beta'=1/2$,
while for part (i) Theorem~\ref{main} only says that $\sup_x |F_n-G_n|=\Lo(n^{-j/2})$ which in fact is $\LO(n^{-(j/2+\delta)})$, and probably,
under slightly stronger assumptions, $\LO(n^{-(j+1)/2})$.
So presumably---in view of Table~\ref{tabelrisk},
\begin{equation}
0\leq {\rm exMSE}(c_{j,n})-{\rm exMSE}(c_{{\rm \SSs ex};n}))=\LO(n^{-j-1}),\qquad j=0,1,2
\end{equation}
\begin{Rem}\rm\small
We even conjecture that we may apply an analogue to Proposition~\ref{ImplMinThm} for
functions $F,F_n,G_n\colon \Psi \to \R$: Let us denote by $\hat \psi^{(j;n)}$, the corresponding
f-o, s-o, t-o optimal IC and $\hat \psi^{({\rm \SSs ex};n)}$ the  exactly optimal IC; then,  with the usual abuse of notation as to ${\rm exMSE}$,
we conjecture that
\begin{equation}
0\leq {\rm exMSE}(\hat \psi^{(j;n)})-{\rm exMSE}(\hat \psi^{({\rm \SSs ex};n)})=\LO(n^{-j-1}),\qquad j=0,1,2
\end{equation}
\end{Rem}
\subsection{Relative risk}\label{relrisk}
An observation in the simulation study was that the relative MSE w.r.t.\ the MSE of the f-o-o procedure
seemed to converge faster than the absolute terms.
This is reflected by our formulas as follows:
\subsubsection{Contaminated situation}
\begin{figure}
  \begin{center}
    \includegraphics[width=10cm,height=7.cm]{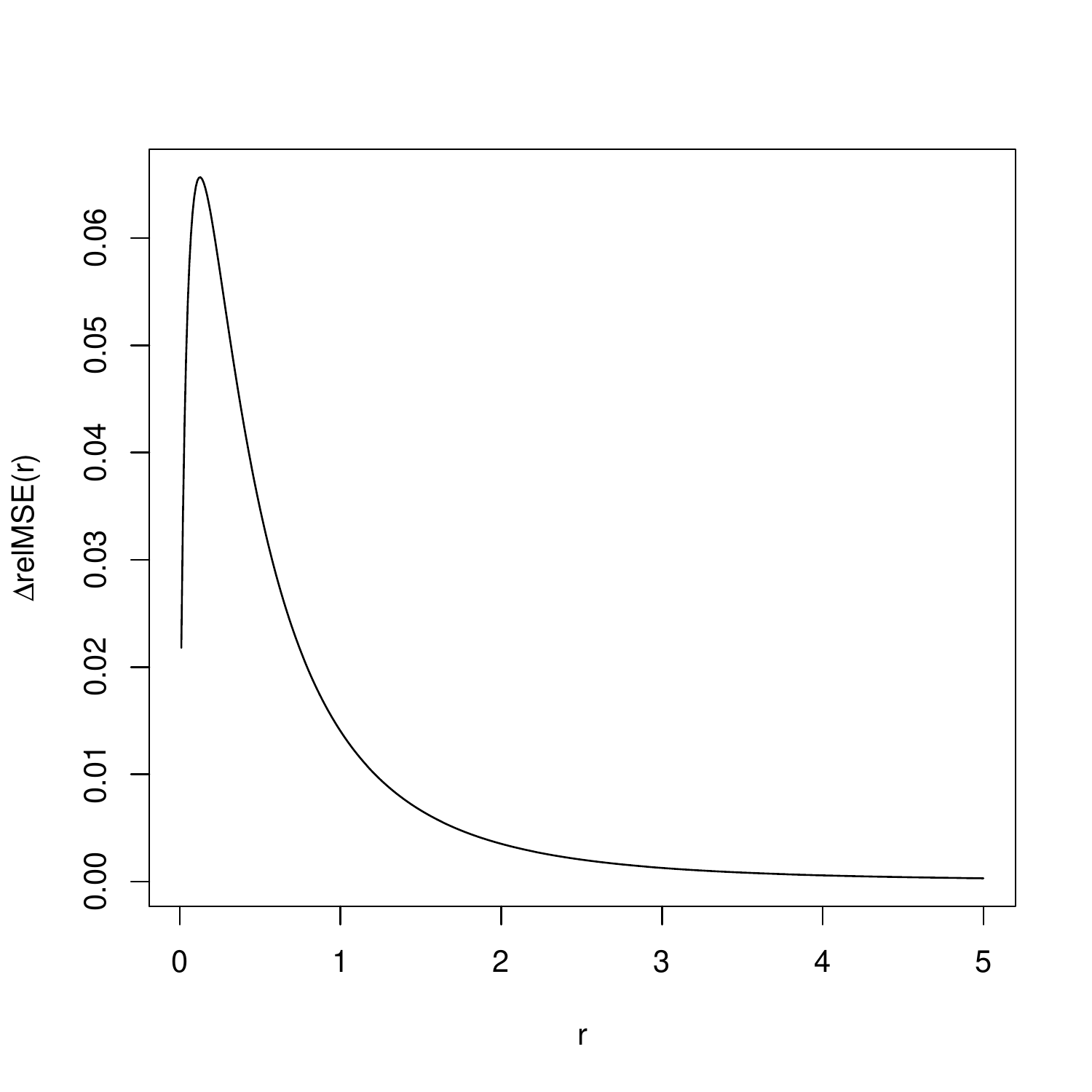}
    \caption{\label{fig1co}{\rm\small The mapping $r\mapsto \widehat{\Delta{\rm relMSE}}(r;\rho)$
    for $F={\cal N}(0,1)$ and for $\rho=0.1$.}}
  \end{center}
\end{figure}

Let ${\rm asMSE}_0(c)$ and $A_1(c)$ be the f-o \asy MSE and the corresponding s-o correction term for the
Hampel-IC with clipping height $c$. Then
we may write for the f-o [s-o]  relative risk ${\rm relMSE}_0(c,r)$ [${\rm relMSE}_1(c,r,n)$] w.r.t.\ the corresponding
risk of the f-o-o procedure
\begin{eqnarray}
{\rm relMSE}_1(c,r,n)&:=&\frac{{\rm asMSE}_0(c)+\Tfrac{r}{\sqrt{n}\,}{A}_1(c)}{{\rm asMSE}_0(c_0)+\Tfrac{r}{\sqrt{n}\,}{A}_1(c_0)}=\nonumber\\
&=&{\rm relMSE}_0(c,r)\left(1+\frac{r}{\sqrt{n}\,}(\Delta(c)-\Delta(c_0))
\right)+\Lo(n^{-1/2})
\end{eqnarray}
with \begin{equation}
  \Delta(c):=\frac{b^2(c)-v_0^2(c)}{{\rm asMSE}_0(c)}
\end{equation}
So in fact, the observed faster convergence is not reflected by higher order optimality, but
as we will see, the difference between ${\rm relMSE}_0(c,r)$ and
${\rm relMSE}_1(c,r)$ are in fact small.\\
Procedure choice will usually be based on relative risk, so it is interesting to consider the maximal
error compared to the s-o approximation one incurs when using the f-o asymptotics instead.
In view of subsection~\ref{soao} we will limit ourselves to only considering Hampel-IC's with
a clipping height $c$ in the range
\begin{equation}
C(c_0,\rho):=[c_0/(1+\rho),\,c_0(1+ \rho)],
\end{equation} for $\rho \geq 0$. This leads us to
\begin{equation}
\widehat{\Delta{\rm relMSE}}(r;\rho):=\max_{c\in C(c_0(r),\,\rho)}r\,\, \big(\Delta(c)-\Delta(c_0(r))\big)
\end{equation}
or even maximizing over the radius
\begin{equation}
\widehat\Delta (\rho):=\widehat{\widehat{\Delta{\rm relMSE}}}(\rho):=\max_{r}\widehat{\Delta{\rm relMSE}}(r;\rho)
\end{equation}
In the Gaussian case, the function $r\mapsto \widehat{\Delta{\rm relMSE}}(r;\rho)$
is plotted for $\rho=0.1$ in Figure~\ref{fig1co}, 
and for $\widehat\Delta(0.1)$,
we get a value of $0.065$,
which for an actual sample size $n$ has to be divided by $\sqrt{n}$---an astonishingly good approximation!\\
{\bf So down to very moderate sample sizes we can base our decision which clipping height to take to achieve
``nearly'' the optimal MSE on $\tilde {\cal Q}_n$ on f-o asymptotics only.} 
%
A similar consideration is of course possible for the ideal situation.

%
%
\subsubsection{Illustration}\label{illussec}
\begin{figure}
  \begin{center}
    \includegraphics[width=12cm,height=7.5cm]{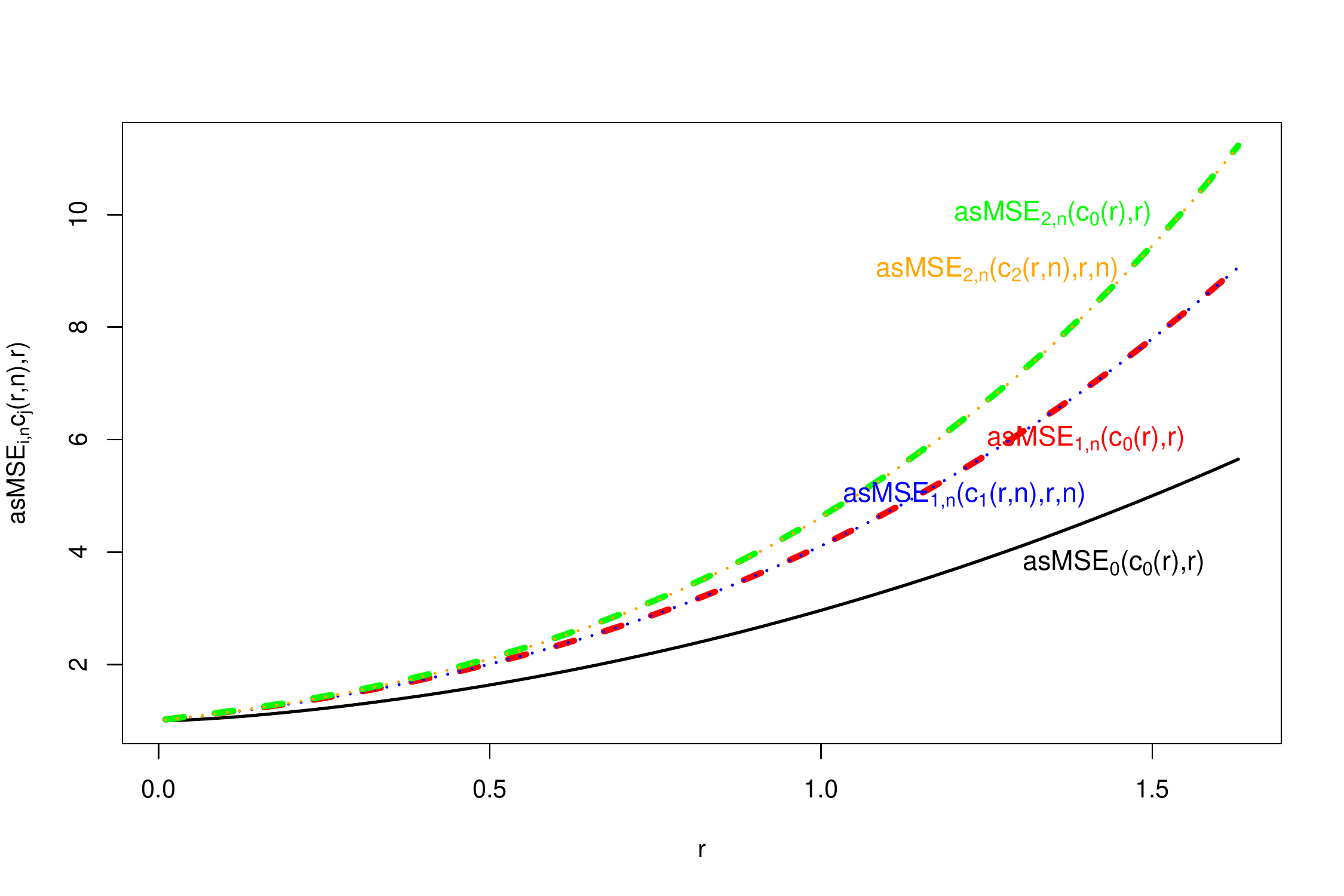}
    \parbox[t]{10cm}{
    \caption[The mapping $r\mapsto {\rm asMSE}$]%
    {\label{figm2}{\parbox[t]{8cm}{{\rm\small The mapping $r\mapsto {\rm asMSE}_{i[,n]}(\eta_{c_j(r[,n])},r[,n])$
    for $i=0,1,2$, $j=0, i$, $n=30$ and $F={\cal N}(0,1)$}}}}}
  \end{center}
\end{figure}

\begin{figure}
  \begin{center}
    \includegraphics[width=10cm,height=7.5cm]{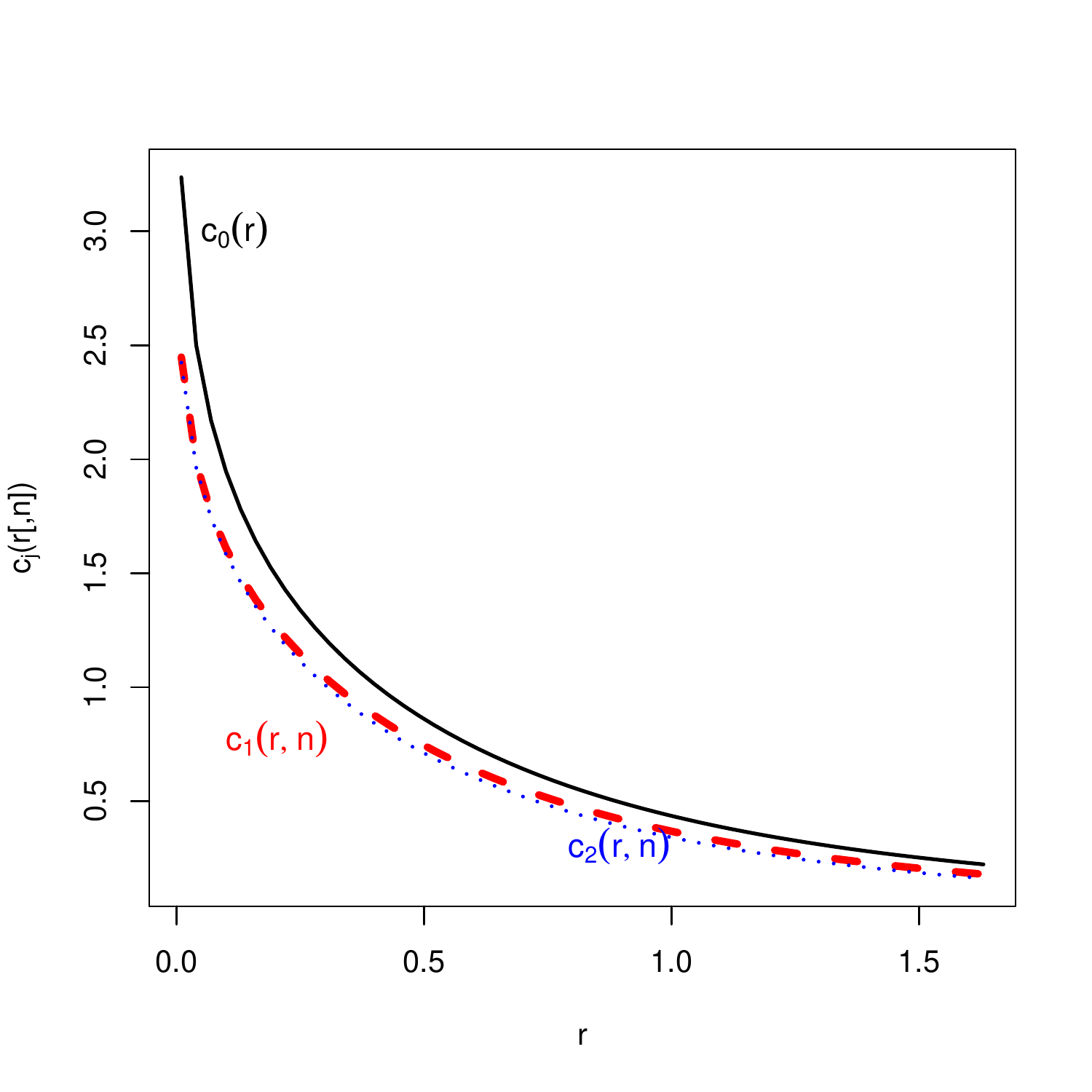}
    \parbox[t]{10cm}{
    \caption[The mapping $r\mapsto c_j$]%
    {\label{figr2}{\parbox[t]{8cm}{{\rm\small The mapping $r\mapsto c_j(r[,n])$
    for $j=0,1,2$ and $n=30$ and $F={\cal N}(0,1)$}}}}}
  \end{center}
\end{figure}

As an example we take $F={\cal N}(0,1)$ and calculate the terms $c_1$,
\begin{equation}
{\rm asMSE}_1:={\rm asMSE}_0+\Tfrac{r}{\sqrt{n}}A_1
\end{equation} and
${\rm relMSE}_1$ for the radii and sample sizes of the simulation study
where for the optimization for $c_1$ we use the function {\tt optimize} in {\tt R 2.11.0}
(compare \citet{RMANUAL}).
The results are tabulated in Table~\ref{tabel12}. 
Correspondingly, we also determine the t-o terms $c_2$,
\begin{equation}
{\rm asMSE}_2:={\rm asMSE}_1+A_2/n
\end{equation}
 and in Figure~\ref{figm2},
we plot the graphs of the five functions
\begin{align*}
r&\mapsto {\rm asMSE}_{0}(\eta_{c_0(r)},r),&
r&\mapsto {\rm asMSE}_{1}(\eta_{c_0(r)},r,n),&
r&\mapsto {\rm asMSE}_{2}(\eta_{c_0(r)},r,n)\\
r&\mapsto {\rm asMSE}_{1}(\eta_{c_1(r,n)},r,n),&
r&\mapsto {\rm asMSE}_{2}(\eta_{c_2(r,n)},r,n)
\end{align*}
for $F={\cal N}(0,1)$ and for 
$n=30$. 
 In fact, the choice of the clipping height---$c_0(r)$, $c_1(r,n)$, $c_2(r,n)$---does not entail
 any visible changes while the absolute value of f-o, s-o, and t-o MSE clearly differ.\\
In the same situation,
the three functions
$r\mapsto c_0(r)$,
$r\mapsto c_1(r,n)$, 
$r\mapsto c_2(r,n)$ 
are plotted in Figure~\ref{figr2}; 
while there are visible differences between $c_0(r)$
and  $c_{i}(r,n)$, $i=1,2$, $c_{1}(r,n)$ and $c_2(r,n)$ visually coincide.
\begin{table}
\caption{\label{tabel12}$c_1(r,n)$, ${\rm asMSE}_1(c_1(r,n),r,n)$ and ${\rm relMSE}_1(c_1(r,n),r,n)$}
\begin{center}
\begin{tabular}{ll||c|c|c|c|c|c|}
$r$ & &
$n=5$&
$n=10$&
$n=30$&
$n=50$&
$n=100$&
$n=\infty$\\
\hline
         &$c_1$           & 1.394  & 1.484  & 1.611  & 1.663  & 1.724  & 1.948\\
{$  0.1$}&${\rm asMSE}_1$ & 1.248  & 1.197  & 1.140  & 1.122  & 1.103  & 1.054\\
         &${\rm relMSE}_1$& 3.476\%& 2.149\%& 0.939\%& 0.623\%& 0.349\%& 0.000\%\\
\hline
         &$c_1$           & 0.994  & 1.059  & 1.147  & 1.181  & 1.219  & 1.339\\
{$ 0.25$}&${\rm asMSE}_1$ & 1.635  & 1.519  & 1.397  & 1.358  & 1.319  & 1.220\\
         &${\rm relMSE}_1$& 2.377\%& 1.470\%& 0.632\%& 0.414\%& 0.228\%& 0.000\%\\
\hline
         &$c_1$           & 0.650  & 0.690  & 0.746  & 0.767  & 0.790  & 0.862\\
{$ 0.5$} &${\rm asMSE}_1$ & 2.527  & 2.271  & 2.006  & 1.923  & 1.840  & 1.636\\
         &${\rm relMSE}_1$& 1.214\%& 0.772\%& 0.342\%& 0.226\%& 0.126\%& 0.000\%\\
\hline
         &$c_1$           & 0.320  & 0.340  & 0.369  & 0.380  & 0.394  & 0.436\\
{$ 1.0$} &${\rm asMSE}_1$ & 5.761  & 4.944  & 4.110  & 3.852  & 3.593  & 2.964\\
         &${\rm relMSE}_1$& 0.427\%& 0.292\%& 0.142\%& 0.098\%& 0.056\%& 0.000\%\\
\end{tabular}
\end{center}
\end{table}
%
%
%
%
\subsection{Comparison with the approach by Fraiman et al.\ (2001)}
\citet{F:Y:Z:01} work in a similar setup, i.e.\ the one-dimensional location problem where the center distribution
is $F_0={\cal N}(0,\sigma^2)$ and an M-estimator $S_n$ to skew symmetric scores $\psi$ is searched which minimizes the maximal risk
on a neighborhood about $F_0$.
Contrary to our approach, the authors work with convex contamination neighborhoods ${\cal V}={\cal V}(F,\ve)$ to a fixed radius $\ve$.\\
There has been some discussion which approach---fixed or shrinking radius---is more appropriate, but for fixed sample size $n$, of course
we may translate the fixed radius $\ve$ into our radius $r/\sqrt{n}$ and then compare the approximation quality of both approaches.\\
\citet{F:Y:Z:01} propose to use risks which are constructed by means of a positive function $g:\R\times \R_+\to\R_+$ of \asy
bias $B=B(F,\psi)$  and \asy variance $v^2=V^2(F,\psi)$. Here,  $B$ is defined as zero of $\beta\mapsto (1-\ve) \int\psi_{\beta}\, dF+\ve b$,
and 
  $v^2:=V_1/V_2^2$ for $V_1=(1-\ve)\int \psi_{B}^2\, dF+\ve b^2$ and $V_2=(1-\ve) \int \dot\psi_{B}\, dF$.

Function $g$ is assumed lower semicontinuous and symmetric in the first argument as well as isotone in each argument.
The risk of an M-estimator to IC $\psi$ is taken as the function
\begin{equation}\label{FYZr1}
  L_g(\psi)=\sup_{G\in {\cal V}} g(B(G,\psi),v(G,\psi)/n)
\end{equation}
A MSE-type risk then is given by $g(u,v)=u^2+v$. It is not quite MSE, as it employs the \asy terms
$B$ and $v$, so their results may differ from ours. The crucial point is that to solve their optimization problem, the authors
have to assume that besides bias, also variance is maximized (for their optimal $\hat\psi$) if we contaminate
with a Dirac measure in $\infty$. According to this assumption, if we introduce $G_0:=(1-\ve)F_0+\ve\Jc_{\{\infty\}}$,
we have to find $\psi$ minimizing
\begin{equation}\label{FYZr2}
  l_g(\psi)=g(B(G_0,\psi),v(G_0,\psi)/n)
\end{equation}
Differently to the Hampel-type IC's the solutions to this problem are of form
\begin{eqnarray}
\psi_{a,b,c,t}(x)&=&\tilde \psi_{a,b,t}\big(x \min\{1,\Tfrac{c}{|x|}\}\big),\\
 \tilde \psi_{a,b,t}(x)&=&a \tanh(tx) + b [x-t\tanh(tx)]\label{FYZs}
\end{eqnarray}
but the ``MSE''-optimal solutions are numerically quite close to corresponding Hampel-ICs $\psi_H$, for which the authors in turn
show that always $L_g(\psi_H)=l_g(\psi_H)$.\\
For an implementation of this optimization see the {\tt R}-file {\tt FYZ.R} available on the web-page.
\rm\smallskip\\
%
\noindent {\bf A comparison}\rm\smallskip\\
As a sort of benchmark for our results, we reproduce a comparison to be found in \citet{Ru:Ko:04}---albeit in some more detail
than in the cited reference:
For a set of values for $n$ and $r$, we determine the ``MSE''-optimal
$\hat\psi$ and a corresponding Hampel IC $\hat \psi_H$ which is then compared to the f-o-o and s-o-o IC derived in this paper.
Within the class of Hampel-IC's, numerically,  we also determine the t-o-o and the ``exactly'' optimal clipping-$c$, $c_2$ and
$c_{\SSs \rm ex}$ respectively. We compare the resulting IC's as to their clipping-height and the corresponding (numerically exact)
value of $R_n(S_n,r)$, denoted by ${\rm MSE}_n$; the latter comparison is done by the terms ${\rm relMSE}^{\SSs \rm ex}_n(c_{\cdot})$, calculated
as \begin{equation}
{\rm relMSE}^{\SSs \rm ex}_n(c_{\cdot}) = ({{\rm MSE}_n(c_{\cdot})}\Big/{{\rm MSE}_n(c_{\SSs \rm ex})}-1)\times 100\%
\end{equation}
The results are displayed in Table~\ref{tabelrisk}. Also compare the function {\tt allMSEs} in the {\tt R}-file {\tt asMSE.R} available on the web-page
to this article.\\
For the numerical evaluation of the ${\rm MSE}$, we use Algorithms C (more accurate, but slow for larger $n$) and D (a little inaccurate for
small $n$, but fast) discussed in \citet{Ru:Ko:04}.
For $n=\infty$, we evaluate the corresponding f-o \asy ${\rm MSE}$ for the IC to the corresponding values of $c$.
As a cross-check, the clipping heights $c_i$, $i=0,1,2$ are also determined for $n=10^8$.
In case of $c_{\SSs \rm FZY}$, for all finite $n$'s the error tolerance used in {\tt optimize} in {\tt R} was $10^{-4}$,
while for $n=\infty$ it was $10^{-12}$. For $c_{\SSs \rm ex}$ and $n=10^8$, an optimization of the (numerically) exact ${\rm MSE}$
would have been too time-consuming and has been skipped for this reason. Also, for $n=5$, the radius $r=1.0$, corresponding to $\ve=0.447$,
is not admitted for an optimization of \eqref{FYZr2} and thus no result is available in this case.
\begin{table}[p]
\caption{\label{tabelrisk}Optimal clipping heights and corresponding (numerically) exact ${\rm MSE}$}
\begin{small}
\begin{center}
\begin{tabular}{ll||c|c|c|c|c|c|}
$r$ & &
$n=5$&
$n=10$&
$n=30$&
$n=50$&
$n=100$&
$n=\infty$\\
\hline\hline
         &$c_0$                             &   1.948  & 1.948  & 1.948  & 1.948  & 1.948  & 1.948   \\
         &${\rm relMSE}^{\SSs \rm ex}_n(c_0)$             & 8.679\% &   4.065\% &   1.340\% &   0.836\% &   0.448\%&{\rm --}\\
\cline{2-8}
         &$c_1$                             &   1.394  & 1.484  & 1.611  & 1.663  & 1.724  & 1.948   \\
         &${\rm relMSE}^{\SSs \rm ex}_n(c_1)$             & 0.833\% &   0.207\% &   0.027\% &   0.014\% &   0.010\%&{\rm --}\\
\cline{2-8}
\raisebox{-1.5ex}[1.5ex]{$  0.1$}&$c_2$                             &   1.309  & 1.428  & 1.585  & 1.644  & 1.713  & 1.948   \\
         &${\rm relMSE}^{\SSs \rm ex}_n(c_2)$             & 0.332\% &   0.066\% &   0.008\% &   0.004\% &   0.006\%&{\rm --}\\
\cline{2-8}
         &$c_{\SSs \rm FZY}$                &   1.368  & 1.370  & 1.610  & 1.668  & 1.756  & 1.939   \\
         &${\rm relMSE}^{\SSs \rm ex}_n(c_{\SSs \rm FZY})$             & 0.658\% &   0.002\% &   0.026\% &   0.021\% &   0.031\%&{\rm --}\\
\cline{2-8}
         &$c_{\SSs \rm ex}$                 &   1.167  & 1.358  & 1.560  & 1.630  & 1.704  & {\rm --}\\
         &${\rm MSE}_n(c_{\SSs \rm ex})$ &   1.388  & 1.239  & 1.151  & 1.129  & 1.107  & {\rm --}\\
\hline\hline
         &$c_0$                             &   1.339  & 1.339  & 1.339  & 1.339  & 1.339  & 1.339   \\
         &${\rm relMSE}^{\SSs \rm ex}_n(c_0)$             & 6.280\% &   3.681\% &   1.108\% &   0.656\% &   0.330\%&{\rm --}\\
\cline{2-8}
         &$c_1$                             &   0.994  & 1.059  & 1.147  & 1.181  & 1.219  & 1.339   \\
         &${\rm relMSE}^{\SSs \rm ex}_n(c_1)$             & 0.933\% &   0.415\% &   0.055\% &   0.023\% &   0.009\%&{\rm --}\\
\cline{2-8}
\raisebox{-1.5ex}[1.5ex]{$ 0.25$}&$c_2$                             &   0.890  & 0.990  & 1.114  & 1.159  & 1.207  & 1.339   \\
         &${\rm relMSE}^{\SSs \rm ex}_n(c_2)$             & 0.241\% &   0.104\% &   0.009\% &   0.002\% &   0.003\%&{\rm --}\\
\cline{2-8}
         &$c_{\SSs \rm FZY}$                &   0.924  & 1.020  & 1.205  & 1.177  & 1.211  & 1.338   \\
         &${\rm relMSE}^{\SSs \rm ex}_n(c_{\SSs \rm FZY})$             & 0.417\% &   0.215\% &   0.233\% &   0.018\% &   0.002\%&{\rm --}\\
\cline{2-8}
         &$c_{\SSs \rm ex}$                 &   0.783  & 0.921  & 1.092  & 1.140  & 1.205  & {\rm --}\\
         &${\rm MSE}_n(c_{\SSs \rm ex})$ &   2.225  & 1.705  & 1.438  & 1.381  & 1.330  & {\rm --}\\
\hline\hline
         &$c_0$                             &   0.862  & 0.862  & 0.862  & 0.862  & 0.862  & 0.862   \\
         &${\rm relMSE}^{\SSs \rm ex}_n(c_0)$             & 2.930\% &   2.655\% &   0.792\% &   0.446\% &   0.218\%&{\rm --}\\
\cline{2-8}
         &$c_1$                             &   0.650  & 0.690  & 0.746  & 0.767  & 0.790  & 0.862   \\
         &${\rm relMSE}^{\SSs \rm ex}_n(c_1)$             & 0.756\% &   0.615\% &   0.087\% &   0.036\% &   0.013\%&{\rm --}\\
\cline{2-8}
\raisebox{-1.5ex}[1.5ex]{$  0.5$}&$c_2$                             &   0.547  & 0.620  & 0.712  & 0.744  & 0.777  & 0.862   \\
         &${\rm relMSE}^{\SSs \rm ex}_n(c_2)$             & 0.230\% &   0.191\% &   0.015\% &   0.008\% &   0.003\%&{\rm --}\\
\cline{2-8}
         &$c_{\SSs \rm FZY}$                &   0.539  & 0.632  & 0.716  & 0.749  & 0.782  & 0.866   \\
         &${\rm relMSE}^{\SSs \rm ex}_n(c_{\SSs \rm FZY})$             & 0.200\% &   0.248\% &   0.021\% &   0.011\% &   0.008\%&{\rm --}\\
\cline{2-8}
         &$c_{\SSs \rm ex}$                 &   0.413  & 0.531  & 0.686  & 0.728  & 0.770  & {\rm --}\\
         &${\rm MSE}_n(c_{\SSs \rm ex})$ &   4.632  & 3.039  & 2.162  & 2.008  & 1.879  & {\rm --}\\
\hline\hline
         &$c_0$                             &   0.436  & 0.436  & 0.436  & 0.436  & 0.436  & 0.436   \\
         &${\rm relMSE}^{\SSs \rm ex}_n(c_0)$             & 2.716\% &   3.132\% &   0.746\% &   0.348\% &   0.149\%&{\rm --}\\
\cline{2-8}
         &$c_1$                             &   0.320  & 0.340  & 0.369  & 0.380  & 0.394  & 0.436   \\
         &${\rm relMSE}^{\SSs \rm ex}_n(c_1)$             & 1.411\% &   1.610\% &   0.251\% &   0.076\% &   0.021\%&{\rm --}\\
\cline{2-8}
\raisebox{-1.5ex}[1.5ex]{$  1.0$}&$c_2$                             &   0.255  & 0.291  & 0.342  & 0.361  & 0.382  & 0.436   \\
         &${\rm relMSE}^{\SSs \rm ex}_n(c_2)$             & 0.876\% &   0.999\% &   0.123\% &   0.027\% &   0.006\%&{\rm --}\\
\cline{2-8}
         &$c_{\SSs \rm FZY}$                & {\rm --} & 0.281  & 0.344  & 0.375  & 0.387  & 0.440   \\
         &${\rm relMSE}^{\SSs \rm ex}_n(c_{\SSs \rm FZY})$             &    {\rm --} &   0.892\% &   0.132\% &   0.063\% &   0.012\%&{\rm --}\\
\cline{2-8}
         &$c_{\SSs \rm ex}$                 &   0.001  & 0.125  & 0.286  & 0.334  & 0.366  & {\rm --}\\
         &${\rm MSE}_n(c_{\SSs \rm ex})$ &  12.627  & 8.445  & 4.948  & 4.296  & 3.787  & {\rm --}\\
\hline\hline
\end{tabular}
\end{center}
\begin{footnotesize}
\begin{tabular}{l|l|l|l}
$c$&order&determined by&optimal among M-estimators\\
\hline
$c_0$ &f-o-o& num. solution of \eqref{HK3}& to any IC\\
$c_1$ &s-o-o& num. solution of \eqref{saoc}& in ${\cal S}_2$ (see section~\ref{soao})\\
$c_2$ &t-o-o& num. optimization of \eqref{mainres} & in ${\cal H}$ (see section~\ref{soao})\\
$c_{\SSs \rm FZY}$ &---& num. optimization of \eqref{FYZr2} & to \eqref{FYZs}-type ICs\\
$c_{\SSs \rm ex}$ &---& num. optimization of the (num.) exact MSE& in ${\cal H}$ (see section~\ref{soao})
\end{tabular}
where 
\eqref{saoc} is
the s-o analogue to \eqref{HK3}, which is derived in Corollary~\ref{cneqcor}.
\end{footnotesize}
{\footnotesize A more detailed description to this table is located on page~\pageref{FYZr2}.}
\end{small}
\end{table}
\makeatletter
\long\def\@makecaption#1#2{%
  \vskip\abovecaptionskip
  \sbox\@tempboxa{#1: #2}%
  \ifdim \wd\@tempboxa >\hsize
    #1: #2\par
  \else
    \global \@minipagefalse
    \hb@xt@\hsize{\box\@tempboxa\hfil}%
  \fi
  \vskip\belowcaptionskip}
\makeatother

%
%
%
%
\section{Ramifications: Minimax radius and Cniper contamination} \label{ramifsec}
\subsection{Minimax radius}\label{minmaxradiussec}
%
In this subsection, we refine the results of \citet{R:K:R:08}. In the cited paper, we want to give
a guideline to the statistician which procedure to choose if he knows that there is contamination but does not know the radius exactly:
To this end, we consider the maximal inefficiency  $\bar \rho(r')$ defined as
\begin{equation} \label{rhobardef}
\bar \rho_0(r'):=  \sup_{r\in(r_l,r_u)} \bar \rho(r',r),\qquad  \bar \rho(r',r):= \frac{\bar R(\eta_{c_0(r')},r)}{\bar R(\eta_{c_0(r)},r)}
\end{equation}
and determine the minimax radius $r_0$ as minimizer of $\bar\rho_0(r')$. If one knows at least that the actual radius will lie in an
interval $[r/\gamma,r\gamma]$ we may determine $r_{\gamma,r}$ as  minimizer of $\bar\rho_\gamma(r',r)=\sup_{s\in(r/\gamma,r\gamma)} \bar\rho(r',s)$ and denote
the corresponding minimax inefficiency by $\bar \rho_{\gamma}(r)$. In a second optimizing step we then determine the maximizer $r_{\gamma}$ of
$\bar \rho_{\gamma}(r)$. The unrestricted case is symbolically included by $\gamma=\infty$. In the Gaussian location case this gives \rm\smallskip\newline
\centerline{\begin{tabular}{||c|c|c||c|c|c||c|c|c||}
\multicolumn{3}{||l||}{$\gamma=0$}&\multicolumn{3}{l||}{$\gamma=2$}&\multicolumn{3}{l||}{$\gamma=3$}\\
$r_0$&$c_0(r_0)$&$\bar \rho_0(r_0)$&$r_2$&$c_0(r_2)$&$\bar \rho_2(r_2)$&$r_3$&$c_0(r_3)$&$\bar \rho_3(r_3)$\\
\hline
$0.621$&$0.718$&$18.07\%$&$0.575$&$0.769$&$8.84\%$&$0.549$&$0.799$&$4.41\%$
\end{tabular}}\medskip\newline
%
%
These calculations can easily be translated to the s-o setup
setting
\begin{equation} \label{rhotildedef}
R_1(\psi,r,n):=r^2 \sup |\psi|^2+ \Ew \psi^2 + \Tfrac{r}{\sqrt{n}}A_1
\end{equation}
so that in this paper we would instead determine $r_1(n)$ as minimizer of $\rho_1(r',r,n)$,
\begin{equation}
  \sup_{r\in(r_l,r_u)}\rho_1(r',r,n),\qquad  \rho_1(r',r,n):=\frac{R_1(\eta_{c_1(r'(n),n)},r,n)}{R_1(\eta_{c_1(r,n)},r,n)}
\end{equation}
respectively $ \rho_{1;\gamma}$ and 
instead of $\bar \rho_{\gamma}$. For finite $n$, however, we have to take into account that
$r<\sqrt{n}$ always. Doing so we get Table~\ref{tabel13}, showing that there is not much variation in both
$c_1(r_{\infty},\cdot)$, $ \rho_{1;\gamma}(r_{\gamma},\cdot)$ for varying $n$.\\
{\bf So if $r$ is completely unknown, it is a good choice to
use the M-estimator to Hampel-scores for $c\approx 0.7$---you will never have a larger inefficiency than the limiting $18\%$!}
Ex post this is one more argument, why the {\tt H07}-estimate survived in in Sections~7.B.8 and 7.C.4 of the Princeton robustness
study (\citet{A:B:H:H:R:T:72}). 
A table for the corresponding t-o minimax radii is available on the web-page.
\begin{table}
\caption{\label{tabel13} Minimax radii for second order asymptotics}%
\begin{center}
\begin{tabular}{ll||c|c|c|c|c|c|}
& &
$n=5$&
$n=10$&
$n=30$&
$n=50$&
$n=100$&
$n=\infty$\\
\hline
               &$r_{\gamma}$                & 0.390  &  0.449  &   0.514  &   0.536  &   0.559  &   0.621\\
$\gamma=0$&$c_1(r_{\gamma})$           & 0.776  &  0.749  &   0.729  &   0.725  &   0.722  &   0.718\\
         &$ \rho_{1;\gamma}(r_{\gamma})$& 16.27\%&  17.08\%&   17.71\%&   17.85\%&   17.96\%&   18.07\%\\
\hline
               &$r_{\gamma}$                & 0.481  &   0.496  &   0.518  &   0.524  &    0.534  & 0.548  \\
$\gamma=3$     &$c_1(r_{\gamma})$           & 0.670  &   0.694  &   0.724  &   0.739  &    0.750  & 0.800  \\
         &$ \rho_{1;\gamma}(r_{\gamma})$& 6.213\%&   6.773\%&   7.490\%&   7.751\%&    8.036\%& 8.836\%\\
\hline
               &$r_{\gamma}$                &  0.540  &  0.552  &  0.564  &  0.563  &  0.571  &  0.574  \\
$\gamma=2$     &$c_1(r_{\gamma})$           &  0.609  &  0.637  &  0.675  &  0.695  &  0.707  &  0.770  \\
         &$ \rho_{1;\gamma}(r_{\gamma})$&  2.987\%&  3.297\%&  3.692\%&  3.834\%&  3.988\%&  4.410\%
\end{tabular}
\end{center}
\end{table}
%
%
\subsection{Innocent-looking risk-maximizing contaminations}\label{milfsec}
In \citet[p.~62]{Hu:97b}, the author complains ``\ldots the considerable confusion between the respective roles of diagnostics and robustness. The purpose
of robustness is to safeguard against deviations from the assumptions, in particular against those that are near or below the limits of detectability.'' 
As worked out in \citet{Ruck:03e}, the exact critical rate for these limits may be determined in a statistical way: 
For some prescribed outlier set ${\rm OUT}$, let $p_0$ and $q_n=(1-r_n)p_0+r_n$ be the probability under the ideal model,
and under convex contaminations of radius $r_n$, respectively. 
Considering the minimax test between these alternatives yields the exact critical rate  $1/\!\sqrt{n}$: under a faster shrinking $p_0$ cannot be
separated from $q_n$ at all, while at a slower rate, asymptotically we can separate them without error.

Going one step further, for some given $1/\!\sqrt{n}$-shrinking neighborhoods of radius $r$, we would also
like to know how ``small'' an outlier may be, while it is still harmful enough
to distort the classically optimal procedure in a way that this procedure is beaten by some robust one.
%
\subsubsection[The Cniper contaminaton]{The {\it Cniper} contaminaton}
%
%
To a fixed radius $r$, in the preceding sections, we have found/discussed f-o-o and s-o-o ICs of
Hampel-form with clipping height $c_j=c_j(r[,n])$, $j=0,1$. 
To these ICs we have derived families of contaminations achieving maximal risk on $\tilde{\cal Q}_n(r)$.
By means of Theorem~\ref{main}(b), these 
are induced by any contaminating measures $P_n^{\rm \SSs di}$ under which
$\eta_{\theta}(X^{\rm \SSs di})$ is constantly either $b_j$ or $-b_j$ for $b_j=A_jc_j$---up to an event
of probability $\Lo(n^{-1})$.
Out of these risk-maximizing contaminations,
let us limit ourselves to those induced by Dirac masses at $x$:
\begin{equation}
Q_n(x):=[(1-\Tfrac{r}{\sqrt{n}})P_{\theta}+\Tfrac{r}{\sqrt{n}}\Jc_{\{x\}}]^{\otimes\, n}
\end{equation}
Among these 
$Q_n(x)$, we seek the least ``suspicious'' looking contamination point $x$
 in the sense that the region ${\rm OUT}_j:=[x;\infty)$ [or $(-\infty;x)$] carries large
ideal probability. With this region as outlier set in \citet{Ruck:03e},
values of $x$  (or slightly above in absolute value) occurring more frequently than they should under the ideal situation,
are hardest to detect.

More precisely, in the general smooth parametric setup (compare \citet{K:R:R:10}),
assume that the observations are univariate; let $S^{\SSs (b_0)}_n$ and $\hat S_n$ be ALEs to
the classical optimal IC  $\hat\eta={\cal I}^{-1}\Lambda$ and the ${\rm asMSE}_0$-optimal IC $\eta_{b_0}$, respectively.
In this setup we define
\begin{Def}
The {\em f-o cniper point} $x_0$ is defined as $x_{0,+}$ if $x_{0,+}\geq - x_{0,-}$ and $x_{0,-}$ else, where
\begin{equation}
  \begin{array}{rcl}
    x_{0,+}&:=&\inf \{x>0\,\Big|\,{\rm asMSE}_0(S^{\SSs (b_0)}_n,Q_n(x))<{\rm asMSE}_0(\hat S_n,Q_n(x))\}\\[1ex]
    x_{0,-}&:=&\sup \{x<0\,\Big|\,{\rm asMSE}_0(S^{\SSs (b_0)}_n,Q_n(x))<{\rm asMSE}_0(\hat S_n,Q_n(x))\}
  \end{array}
\end{equation}
\end{Def}

\begin{Rem}\rm\small
\begin{ABC}
\item The name {\it cniper\/} point is due to H.~Rieder; it alludes to the fact that this ``Ianus-type'' contamination $Q_n(x_0)$ pretends to
be {\it nice\/},
but to the contrary is in fact {\it pernic}ious, ``sniping'' off the classically optimal procedure\ldots
\item The cniper concept is of course not bound to quadratic loss. 
%
%
In the obvious manor, the concept may be generalized for multivariate observations, if  we define any $x_0$ of minimal absolute as {\it cniper} point.
\item To get rid of the dependency upon the radius $r$, in the examples we will use the minimax radii $r_{\gamma}(n)$ defined in
 the preceding section.
\end{ABC}
\end{Rem}
Correspondingly, in the setup of this paper and under \eqref{symdef},
let $S^{\SSs (c_1)}_n$ be an M-estimator to the s-o-o IC $\eta_{c_1}$ according to Corollary~\ref{cneqcor}.
\begin{Def}
The {\em s-o cniper point} $x_1$ is defined as $x_{1,+}$ if $x_{1,+}\geq - x_{1,-}$ and $x_{1,-}$ else, where
\begin{equation}
  \begin{array}{rcl}
    x_{1,+}&:=&\inf \{x>0\,\Big|\,{\rm asMSE}_1(S^{\SSs (c_1)}_n,Q_n(x))<{\rm asMSE}_1(\hat S_n,Q_n(x))\}\\[1ex]
    x_{1,-}&:=&\sup \{x<0\,\Big|\,{\rm asMSE}_1(S^{\SSs (c_1)}_n,Q_n(x))<{\rm asMSE}_1(\hat S_n,Q_n(x))\}
  \end{array}
\end{equation}
\end{Def}
{\it Cniper} contaminations and f/s-o-o ICs form saddle-points under \eqref{minetahat}/\eqref{symdef}:
\begin{Prop}\label{sattel}
The pair $(S^{\SSs (b_0)}_n,Q_n(x_0))$ is a saddlepoint
for the class of all pairs $(S_n,Q_n)$
 if
\begin{equation}\label{minetahat}
  |\hat \eta(x_0)|\leq |\eta_b(x_0)|\qquad\forall b\colon \quad |\eta_b(x_0)|< b
\end{equation}
where $S_n$ are ALE's to IC's of form~\eqref{allgemoptIC}  and $Q_n\in {\cal Q}_n$ w.r.t.\ f-o risk
$\tilde R$.\\
 Under \eqref{symdef}, the same holds in the one-dimensional location model for the pair $(S^{\SSs (c_1)}_n,Q_n(x_1))$ w.r.t.
s-o risk in $\tilde {\cal Q}(r)$.
\end{Prop}
%
%
%
\begin{Rem}\rm\small
A sufficient condition for \eqref{minetahat} 
is that $\Lambda(x)=-\Lambda(-x)$: Then for any $b>0$, $a_b=0$ is possible and, 
$$A_b^{-1}=\Ew \Lambda \Lambda^\tau \min\{1,\frac{b}{|A_b\Lambda|}\}\preceq \Ew \Lambda\Lambda^\tau= {\cal I}$$
So $A_b\succeq {\cal I}^{-1}$ in
the positive semi-definit sense, and hence
for $b$ s.t.\ $|\eta_b(x_j)|< b$
\begin{equation}
|\eta_b(x_j)|=|A_b \Lambda(x_j)|\ge  |{\cal I}^{-1}\Lambda(x_j)|=|\hat \eta(x_j)|
\end{equation}
\end{Rem}
%
\subsubsection{Error probabilities}
%
For numerical evaluations, we consider
the Gaussian location model and   the Gaussian location and scale model.
In both models, $x_{j,+}=-x_{j,-}$, and without loss, we use $x_{j,+}$.\\
For the \asy tests between
$q_n=p_0$ and $q_n> p_0$, alluded to in the beginning of this section,
we note that
\begin{equation}
p_0=P_{\theta}(X_i\geq x_j)=\Phi(-x_j),\qquad q_n=p_0+\frac{r}{\sqrt{n}}(1-p_0)
\end{equation}
As to the (f-o) \asy minimax test \citet[formula (6.1)]{Ruck:03e} gives
as \asy risk
\begin{equation}
\ve=\ve_{\infty}=\Phi\Big(-\frac{r}{2}\,\sqrt{\frac{1-p_0}{p_0}}\Big)
\end{equation}
For s-o asymptotics, we instead use the finite-sample
minimax test, i.e.\ the Neyman-Pearson test with equal Type-I and Type-II error.
In our case this is a corresponding randomized
binomial test. 
\subsubsection{Gaussian location}
In the Gaussian location model, we draw all necessary expressions from \citet[Prop.~]{Ruck:03c}; in particular, with
$c_1=c_1(n,r_{\gamma})$, and $A_1=(2\Phi(c_1)-1)^{-1}$,
$b_1=c_1A_1$, by Theorem~\ref{main}(b),
 maximizing risk amounts to
either $X^{\rm \SSs di}> c_1$ always or $X^{\rm \SSs di}<-c_1$ always.
The classically optimal estimator is the  arithmetic mean, and one easily calculates
\begin{equation}
  \Ew_{Q_n(x)}[\bar x_n^2\, \Big|\,K=k]=\frac{1}{n^2}[k^2x^2+(n-k)]
\end{equation}
and integrating out $K$ we get directly
\begin{equation}
n\,  \Ew_{Q_n(x)}[\bar x_n^2]=1-\Tfrac{r}{\sqrt{n}}+x^2(r^2+\Tfrac{r}{\sqrt{n}}-\Tfrac{r^2}{n})
\end{equation}
Combining this with formulas \eqref{mainres} and \eqref{A1sym}, for $M_0:={\rm asMSE}_0(S^{\SSs (c_1)}_n)$ we get
\begin{equation}
  x_1^2(n)=\frac{M_0-1+\frac{r}{\sqrt{n}}(M_0+b_1^2(r^2+1)+1)}{r^2(1-\frac{1}{n})+\frac{r}{\sqrt{n}}}
\end{equation}
or
\begin{equation}
  x_1(n)=\frac{\sqrt{M_0-1}}{r}+\frac{1}{2\,\sqrt{n}}[\frac{M_0+1+b_1^2(r^2+1)}{\,\sqrt{M_0-1}\,}-\frac{\sqrt{M_0-1}\,}{r^2}]+\Lo(\Tfrac{1}{\sqrt{n}})
\end{equation}
This yields the results as in Table~\ref{tabel14}. 
We include the type-II error $1-\beta(\alpha)$ for the Neyman Pearson test to niveau $\alpha=5\%$ and the risk $\ve_n$ of the corresponding
minimax test; roughly speaking we cannot do better than overlooking one of $10$
contaminations at niveau $5\%$ ideal observations to be falsely marked as outliers, and, equally weighting the two error types
we cannot do better than with a false classification rate of $7\%$ for each error type.
\begin{table}
\caption{\label{tabel14} Minimax contamination at $\gamma=0$}%
\begin{center}
\begin{tabular}{l||c|c|c|c|c|c|c|c|}
 $n$&
$5$&
$10$&
$30$&
$50$&
$100$&
$200$&
$300$&
$\infty$\\
\hline
$r_{\gamma}(n)$       & 0.390 & 0.449 & 0.514 & 0.536 & 0.559 & 0.576&  0.584 & 0.621\\
$c_1(r_{\gamma},n)$   & 0.776 & 0.749 & 0.729 & 0.725 & 0.722 & 0.720&  0.719 & 0.718\\
$x_1(n)$              & 2.931 & 2.470 & 2.101 & 2.004 & 1.914 & 1.853 & 1.826 & 1.714\\
$1-\beta_n(0.05)$     & 0.364 & 0.272 & 0.215 & 0.183 & 0.162 & 0.133 & 0.132 & 0.101\\
$\ve_n$               & 0.277 & 0.178 & 0.129 & 0.115 & 0.097 & 0.089 & 0.086 & 0.072
\end{tabular}
\end{center}
\end{table}

\subsubsection{Gaussian location and scale}
    To give one more example, consider the one-dimensional location-scale model at central distribution ${\cal N}(0,1)$.
    For this model we have not yet established a s-o \asy theory; for f-o asymptotics, however,
    we may use {\tt R}-programs from the bundle {\tt RobASt}, cf.\ \citet[Appendix~D]{Koh:04di}, and get $r_{\infty}=0.579$,
    \begin{equation}
    \max_{Q_n\in{\cal Q}_n(r_{\infty})}{\rm asMSE}(\eta_{\theta;0},Q_n)=3.123
    \end{equation}
    while ${\cal I}_{\theta}^{-1}\Lambda_{\theta}=(x,\Tfrac{1}{2}(x^2-1))^{\tau}$. This gives $x_0=1.844$---and hence
    $\ve_{\infty}=5.737\%$ and $1-\beta_{\infty}(5\%)=6.557\%$. Condition~\eqref{minetahat} is proved to hold in
    subsection~\ref{minetp}.

\appendix
\section{Proofs}\label{proofsec}\req
%
%
\subsection{A Hoeffding Bound}
%
%
%
\begin{Lem}\label{hoef2}
  Let $\xi_i \iid F$,  $i=1,\ldots,n$ be real--valued random variables, $|\xi_i|\leq 1$
  Then for $\mu=\Ew[\xi_1]$ and $0<\ve<1-\mu$
  \begin{eqnarray}
    P(\frac{1}{n}\sum_i \xi_i -\mu \geq \ve) &\leq& \left\{\left(\frac{\mu}{\mu+\ve}\right)^{\mu+\ve}
    \left(\frac{1-\mu}{1-\mu-\ve}\right)^{1-\mu-\ve}\right\}^n\label{hoe3}
  \end{eqnarray}
\end{Lem}
\begin{proof}{}
  \citet{Hoef:63}, Thm.~1, inequality (2.1).
\hfill\qed\end{proof}
To settle case (II) in the proof of Theorem~\ref{finminmax}, we need the
following sharpening of \citet[Lem.~A.2]{Ruck:03b}
\begin{Lem}\label{binlem}
Let $k_1(n)=1+d_n$ and assume that for some $\delta\in(0,1/4)$,
 \begin{equation} \label{dbeding}
d_n n^{1/4-\delta}\to \infty,\qquad d_n n^{-1/2+\delta} \to 0\qquad\mbox{for }n\to \infty
 \end{equation}
Then if $\liminf_n d_n>0$
 there is some $c>0$ such that
 \begin{equation} \label{crneq1}
\Pr({\rm Bin}(n,r/\sqrt{n}\,)>k_1(n) r\sqrt{n}\,)=\Lo(e^{- c r\sqrt{n}})
 \end{equation}
 and, if $d_n=\Lo(n^0)$, for any $0<\delta_0\leq 2\delta$, it holds that
\begin{equation} \label{crneq2}
\Pr({\rm Bin}(n,r/\sqrt{n}\,)>k_1(n) r\sqrt{n}\,)=
\Lo(e^{-r n^{\delta_0}})
\end{equation}
\end{Lem}
\begin{Rem}\rm\small{}
Even if $d_n$ is increasing at a faster rate than $n^{1/2}$, assertion~\eqref{crneq1} remains true,
as long as $\liminf_n d_n>0$---but this is not needed here.
\end{Rem}
\begin{proof}{}
Let
\begin{equation}\label{calkdef}
{\Ts {\cal K}}_n:=k_1(n)\log k_1(n)+1-k_1(n)=\int_1^{k_1(n)}\log(x)\, dx
\end{equation}
Then ${\Ts {\cal K}}_n>0$, as $\log(x)>0$ for $x>1$ and
By the second assumption in \eqref{dbeding}, $d_n=\Lo(\sqrt{n\,}\,)$, so $0<d_nr/\sqrt{n}<1-r/\sqrt{n}$
eventually in $n$ and Hoeffding's Lemma~\ref{hoef2} is available; applying it
to the case of $n$ independent ${\rm Bin}(1,p)$ variables,
we obtain for $B_n\sim{\rm Bin}(n,p_n)$, $p_n=r/\sqrt{n}$ and
$\ve=(k_1(n)-1)r/\sqrt{n}$ (which is smaller than $1-p_n$ eventually)
\begin{eqnarray*}
  \Pr(B_n>k_1(n)r\sqrt{n}\,)&\leq& 
\exp\Big(-k_1(n)r\sqrt{n}\,\log(k_1(n))+(n-k_1(n)r\sqrt{n}\,)\times\\
 &&\qquad \times\big(\log(1-\frac{r}{\sqrt{n}\,})-\log(1-k_1(n)\frac{r}{\sqrt{n}\,})\big)\Big)
\end{eqnarray*}
But for $x_0<x_1 \in(0,1)$, $\log(1-x_0)-\log(1-x_1)=\int_{1-x_0}^{1-x_1}t^{-1}\,dt\leq (x_1-x_0)/(1-x_1)$.
Thus
$
\log(1-r/\sqrt{n}\,)-\log(1-k_1(n)r/\sqrt{n}\,)\leq \frac{d_n r/\sqrt{n}\,}{1-k_1(n)r/\sqrt{n}\,}
$ and
\begin{eqnarray*}
\Pr(B_n>k_1(n)r\sqrt{n}\,)&\leq&
\exp\Big(-r\sqrt{n\,}\,\big(k_1(n)\,\log(k_1(n))-k_1(n)+1\big)\Big)=
e^{-{\cal K}_n\, r\sqrt{n}},
\end{eqnarray*}
If $\liminf_n d_n>0$, by \eqref{calkdef} $\liminf_n {\Ts {\cal K}}_n>0$, and for any $0<c<\liminf_n {\Ts {\cal K}}_n$,
\eqref{crneq1} follows. If $d_n=\Lo(n^0)$,
we note that
\begin{equation}
{\Ts {\cal K}}_n=(1+d_n)\log(1+d_n)-d_n=d_n^2/2+\Lo(d_n^2)
\end{equation}
which for any $\delta'>0$ entails
$$\Pr({\rm Bin}(n,r/\sqrt{n}\,)>k_1(n) r\sqrt{n}\,)=
\Lo\Big(\exp\big(-\frac{r d_n^2\sqrt{n\,}\,}{2+\delta'}\big)\Big)
$$
Now for $d_n=\Lo(n^0)$, by the first assumption in \eqref{dbeding}, for $0<\delta_0<2\delta$
eventually in $n$, \eqref{crneq2} holds as 
$$
n^{\delta_0}-\frac{d_n^2\sqrt{n\,}\,}{2+\delta'}< n^{2 \delta}(1- \frac{n^{1/2-2\delta}d_n^2}{2+\delta'})\to -\infty
$$
\hfill\qed\end{proof}
Another consequence of the exponential decay of \eqref{crneq1}/\eqref{crneq2} is
that we may neglect values of $K>k_1(n) r\sqrt{n}$  when integrating along $K$.
\begin{Cor} \label{corewk}
Let $K\sim{\rm Bin}(n,r/\sqrt{n}\,)$. Then, in the setup of Lemma~\ref{binlem}, for any $j\in\N$,
\begin{equation}
\Ew[K^j \Jc_{\{X\geq k_1(n) r \sqrt{n}\}}]=\Lo(e^{-rn^{d}})
\end{equation}
for any $0<d<\sqrt{n}$ if $\liminf_n d_n>0$ and any $0<d\leq \delta_0$ if $\lim_n d_n=0$.
\end{Cor}
\begin{proof}{}
$\Ew[K^j \Jc_{\{K\geq k_1(n) r \sqrt{n}\}}]
\leq
n^j\Pr(X>k_1(n)r\sqrt{n})\stackrel{\eqref{crneq1}/\eqref{crneq2}}{=} 
\Lo(e^{-rn^{d}})
$.
\hfill\qed\end{proof}

\subsection{Proof of Theorem~\ref{finminmax}}\label{pfinminmax}
In the risk, we have to treat stochastic arguments in $\Phi$, $\varphi$; this is settled in the following lemma:
\begin{Lem}\label{Philem}
Let $F\colon \R\to\R$ be twice differentiable  with H\"older-continuous second derivative and
$G\colon \R\to\R$ be differentiable  with H\"older-continuous derivative.
Then there is a sequence  $k_1(n)=1+d_n$ with $d_n\to 0$ according to \eqref{dbeding} and some $\eta>0$, such that for all $x,\beta \in \R$
and with $\tilde k=K/\sqrt{n}$,
\begin{equation}
  \Ew[F(x+\beta \tilde k) |K\leq k_1(n) r\sqrt{n}\,]=F(x+\beta r)+F''(x+\beta r)\frac{\beta^2r}{2\sqrt{n}}\,+\Lo(n^{-1/2})
\end{equation}
and
\begin{equation}
  \Ew[G(x+\beta \tilde k) |K\leq k_1(n) r\sqrt{n}\,]=G(x+\beta r)+\LO(n^{-(1+\eta)/4})
\end{equation}
\end{Lem}
\begin{proof}{}

Using the Taylor approximation of $\log(1+x)$, we get for $n$ sufficiently large
\begin{equation}
d_n^2/3\leq d_n^2/2-d_n^3/6\leq {\Ts {\cal K}}_n\leq d_n^2/2
\end{equation}

By  \eqref{crneq2} of Lemma~\ref{binlem},  for some $\delta_0$ and eventually in $n$
we have
$
P(K>k_1(n) r\sqrt{n}\,)\leq \exp(-r n^{\delta_0})
$, and 
by the same argument we also get that
$ 
P(K<(2-k_1(n)) r\sqrt{n}\,)\leq  \exp(-r n^{\delta_0})
$. 
Hence,
\begin{equation}
P(|\tilde k-r|>r d_n)\leq  2 \exp(-r n^{\delta_0})
\end{equation}
Thus, as $F$, $G$ are bounded, the contribution of the set $\{|\tilde k-r|>r d_n\}$ decays exponentially, while
on the complement we have a uniformly bounded Taylor expansion up to order $2$ respectively $1$ for the integrands:
\begin{eqnarray*}
  F(x+\beta \tilde k)&=&F(x+\beta r)+F'(x+\beta r)\beta(\tilde k-r)+F''(x+\beta r)\beta^2(\tilde k-r)^2/2
  +\Lo((\tilde k-r)^{2+\eta})\\
  G(x+\beta \tilde k)&=&G(x+\beta r)+G'(x+\beta r)\beta(\tilde k-r)+\Lo((\tilde k-r)^{1+\eta})
\end{eqnarray*}
Integrating these expansions out in $\tilde k$, we see that the first contribution to the Taylor series
 for $F$ is the quadratic term, which is
$F''(x+\beta r)\frac{\beta^2}{2}\, \Var \tilde k$, and the remainder is $\Lo(n^{-1/2})$. For $G$, the first contribution to the error term is the remainder,
hence of form ${\rm const}|\tilde k-r|^{1+\eta}$. By the H\"older inequality this gives a bound
$ 
{\rm const}\, [\Var \tilde k]^{\frac{1+\eta}{2}}=\LO(n^{-(1+\eta)/4})\vspace{-2ex}
$. 
\hfill\qed\end{proof}
For the proof of Theorem~\ref{finminmax}, we use a tableau like the one of \citet[p.~19]{Ruck:03c}, i.e.,
to derive the result, we partition the integrand according  to 
\label{tabelAUF}
\begin{center}
\begin{tabular}{c||c|c}
&$K< k_1(n)r\sqrt{n}$ &  $k_1(n)r\sqrt{n}\leq K< \ve_0 n $ \\[0.5ex]
\hline\hline
$|t|\leq k_2b^2{\log(n)/n}$& (I) &\\[0.8ex]
\cline{1-2}
$k_2 b^2{\log(n)/n} < |t|$& (III)&  \raisebox{1.5ex}[-1.5ex]{(II)} \\[0.8ex]
\hline
\end{tabular}
\end{center}
with $k_1(n)$ according to \eqref{dbeding}. This time, no integration w.r.t.\  $t$ is needed,  so  case (IV) from \citet{Ruck:03c}
may be canceled, which is why we may dispense of assumption~(Pd) and pass to the unrestricted
neighborhoods ${\cal Q}_n$. Cases (II) and (III) may be taken over unchanged from \citet[Proof of Thm.~3.5]{Ruck:03c},
so we may confine us to case (I):\\
We use $\alpha_1$, $\alpha_2$ from \eqref{alph12} and proceed paralleling the proof in \citet{Ruck:03c}
and get from formula~(A.18) therein that
$\Pr(S_n\leq -\Tfrac{\alpha_1}{\sqrt{n}}\,|\,D_{k,\tilde t}\,)=\tilde G_n(\,-\Tfrac{\alpha_1}{\sqrt{n}}\,)+\LO(n^{-3/2})$. 
So we have to spell out $s_{n,k}(\Tfrac{-\alpha_1}{\sqrt{n}}\,)$, 
which gives 
\begin{equation}
s_{n,k}(\Tfrac{-\alpha_1}{\sqrt{n}}\,)=v_0^{-1}\Big\{(-{t^{\natural}}- \alpha_1)+
\Tfrac{1}{\sqrt{n}}[\Tfrac{\tilde k}{2}\alpha_1
-\alpha_1\tilde v_1 ({t^{\natural}}+\alpha_1)
-\Tfrac{l_2}{2}\alpha_1^2]\Big\}+\Lo(\Tfrac{1}{\sqrt{n}}\,)
\end{equation}
and hence---setting $\tilde s=s_{n,k}(\Tfrac{-\alpha_1}{\sqrt{n}}\,)$ and $\tilde s_1=-(\alpha_1+t^{\natural})/v_0$ as in \citet{Ruck:03c}
\begin{eqnarray}
\Pr(S_n\leq -\Tfrac{\alpha_1}{\sqrt{n}}\,|\,D_{k,\tilde t}\,)&=&
\Phi(\tilde s)-
\varphi(\tilde s)\Tfrac{(\tilde s^2-1)}{6\sqrt{n}}\rho(-\Tfrac{\alpha_1}{\sqrt{n}})
+\Lo(\Tfrac{1}{\sqrt{n}}\,)=\nonumber\\
&=&
\Phi(\tilde s_1)+\Tfrac{\varphi(\tilde s_1)}{2\sqrt{n}v_0}[\alpha_1 \tilde k
-l_2 \alpha_1^2- 2(\alpha_1+t^{\natural})\tilde v_1 \alpha_1
-  v_0\Tfrac{\rho_0}3(\tilde s_{1}^2-1)]
+\Lo(\Tfrac{1}{\sqrt{n}}\,)
\end{eqnarray}
This term is maximized eventually in $n$, if $-t^{\natural}$ is maximal or, essentially equivalent,
all contaminating mass (up to mass $\Lo(n^{-1/2})$) is concentrated left of $\check y_n$ from Section~\ref{HOE}, and then
$
t^{\natural}=k^\natural \check b
$,  
and after the substitution according to $\tilde k:=k/\sqrt{n}$, $k^\natural:=k/\sqrt{\bar n}$,
this  gives with $\tilde s_k=-(\alpha_1+\tilde k \check b )/v_0$
\begin{equation}
\Pr(S_n\leq -\Tfrac{\alpha_1}{\sqrt{n}}\,|\,D_{k,\tilde t=k \check b}\,)=
\Phi(\tilde s_k)+
\Tfrac{\varphi(\tilde s_k)}{2\sqrt{n}v_0}[\alpha_1 \tilde k
-l_2 \alpha_1^2- 2\tilde s_k v_0 \tilde v_1 \alpha_1
-  v_0\Tfrac{\rho_0}3(\tilde s_{k}^2-1)-\tilde k^2 \check b]
+\Lo(\Tfrac{1}{\sqrt{n}}\,)
\end{equation}
Now, by \eqref{s1def64}, it holds that
$s_1=-(\alpha_1+r \check b)/v_0$,
so that by an application of Lemma~\ref{Philem}, for $Q^0_{n;\,-}$ any sequence of measures according to \eqref{contbed1}
\begin{eqnarray*}
Q^0_{n;\,-}(S_n\leq -\Tfrac{\alpha_1}{\sqrt{n}})=\Phi(s_1)+\Lo(\Tfrac{1}{\sqrt{n}}\,)+
\Tfrac{1}{\sqrt{n}}\varphi(s_1)
\big [\Tfrac{r}{2v_0} \alpha_1 -\Tfrac{l_2}{2v_0} \alpha_1^2+  s_1 v_0 \tilde v_1 \alpha_1
-   \Tfrac{\rho_0}6(\tilde s_{1}^2-1)- r \Tfrac{\check b^2}{2v_0^2} s_1 - r^2 \Tfrac{\check b}{2v_0}\big]
\end{eqnarray*}
Correspondingly, we get for any sequence of measures $Q^{+}_n$  according to \eqref{contbed2}
\begin{eqnarray*}
Q^0_{n;\,+}(S_n\geq \Tfrac{\alpha_2}{\sqrt{n}})=\Phi(s_1)+\Lo(\Tfrac{1}{\sqrt{n}}\,)+ \Tfrac{1}{\sqrt{n}}\varphi(s_1)
\big [\Tfrac{r}{2v_0} \alpha_2 +\Tfrac{l_2}{2v_0} \alpha_2^2- s_1 v_0 \tilde v_1 \alpha_2
+   \Tfrac{\rho_0}6(\tilde s_{1}^2-1)- r \Tfrac{\hat b^2}{2v_0^2} s_1 + r^2 \Tfrac{\hat b}{2v_0}\big]
\end{eqnarray*}
We next account for order $\Tfrac{1}{\sqrt{n}}$-terms and get, as $\delta'=\LO(\Tfrac{1}{\sqrt{n}}\,)$
\begin{equation}
Q^0_{n;\,-}(S_n\leq -\Tfrac{\alpha_1'}{\sqrt{n}}\,)=
Q^0_{n;\,-}(S_n\leq -\Tfrac{\alpha_1}{\sqrt{n}}\,) +\delta'\varphi(\Tfrac{a -r\bar b}{v_0})+\Lo(\Tfrac{1}{\sqrt{n}}\,)
\end{equation}
and analogously for $Q^0_{n;\,+}(S_n\geq \Tfrac{\alpha_2'}{\sqrt{n}}\,)$,
so
$\delta'=
\Tfrac{1}{\sqrt{n}}\Big(-\Tfrac{r\delta}{2v_0}-\Tfrac{l_2}{2v_0}(a^2+\delta^2)-\tilde v_1 v_0s_1 \delta-
\Tfrac{\rho_0}{6} (s_1^2-1)+ \Tfrac{r\bar b \delta s_1}{v_0^2}+\Tfrac{r^2\delta}{2v_0}\Big)
%
$ 
and $Q^0_{n;\,-}(S_n\leq -\Tfrac{\alpha_1'}{\sqrt{n}}\,)=Q^0_{n;\,+}(S_n\geq \Tfrac{\alpha_2'}{\sqrt{n}}\,)+\Lo(\Tfrac{1}{\sqrt{n}}\,)$, i.e.,
\begin{eqnarray}
Q^0_{n;\,-}(S_n\leq -\Tfrac{\alpha_1'}{\sqrt{n}}\,)=
\Phi(s_1)+
\varphi(s_1) \Tfrac{1}{\sqrt{n}}\Big[\Tfrac{ra}{2v_0}+2\Tfrac{l_2 a \delta}{v_0}-a s_1\tilde v_1-
 \Tfrac{r(\check b^2+\hat b^2)s_1}{4 v_0^2}+\Tfrac{r^2\bar b}{2v_0}
\Big]+\Lo(\Tfrac{1}{\sqrt{n}}\,)
\end{eqnarray}
\hfill\qed

\subsection{Proof of Corollary~\ref{cneqcor}}\label{cneqseq}
The assumptions of Theorem~\ref{cvthm} are clearly fullfilled. Hence we may start with
the verification \eqref{saoc}:
\begin{eqnarray}
G(w,s)&=&(w^2+s^2)(1+\frac{r}{\sqrt{n}\,}\,)+\frac{r}{\sqrt{n\,}}\,w^2(1+\frac{1}{r^2})\\
\partial_wG(w,s)&=&2w[1+\frac{r}{\sqrt{n}\,}+\frac{r}{\sqrt{n}\,}(1+\frac{1}{r^2})\,], \qquad 
\partial_sG(w,s)=
2s[1+\frac{r}{\sqrt{n}\,}\,]
\end{eqnarray}
and hence, dividing both sides of \eqref{diffarg} by $2\hat A\hat v_0$, we get the assertion.
The LHS of \eqref{saoc} (with or without
factor $1+\frac{r^2+1}{r^2+r\sqrt{n}}$) is isotone, the RHS antitone in $c$. Thus if we insert the factor to
correct the f-o-o clipping height $c_0$ to $c_{1}(n)$, the factor increases the LHS without affecting the RHS. This can only be compensated
for by a decrease of $c_0$ to $c_{1}(n)$. If 
$h(c)$ is differentiable in
$c_0$ with derivative $h'(c_0)$,
\eqref{c1def} is an application of the applying the implicit function theorem:
Let $G(s,c):=r^2 c\,(1+s)-h(c)$. Then $G(0,c_0)=0$. Hence for $s=(r^2+1)/(r^2+r\sqrt{n}\,)$,
up to $\Lo(n^{-1/2})$,
\begin{eqnarray*}
&&  c_{1}(n)+ \Lo(n^{-1/2})=c_0 - \frac{G_s(0,c_0)}{G_c(0,c_0)}s=
  c_0\,\Big(1-\frac{1}{\sqrt{n}}\,\frac{r^3+r}{r^2-h'(c_0)}\,\Big)+\Lo(n^{-1/2})
\end{eqnarray*}
\hfill\qed

\subsection{Proof of Proposition~\ref{ImplMinThm}}\label{ImplMinThmsec}
We apply \citet[Theorem~1.4.7]{Ri:94} to the derivatives; this theorem says that for $\eta\in {\cal C}_1(\R)$
with $\eta(\theta_0)=0$,  $\eta'(\theta_0)\not=0$ for some $\theta_0\in\R$, there exists an open neighborhood $V_0\subset{\cal C}_1(\R)$
such that for every open, connected neighborhood $V\subset V_0$ of $\eta$ there is a unique, continuous map $T\colon V\to \R$ with
\begin{equation} \label{Tabbdiff}
T(\eta)=\theta_0, \qquad f(T(f))=0,\quad f\in V
\end{equation}
even more so, $T$ is continuously bounded differentiable on $V$ with derivative at tangent $h$
\begin{equation} \label{Tabdiff}
dT(f)h=-h(T(f))/f'(T(f))
\end{equation}
Hence there is an open neighborhood $V_{0; F}$ of $F$ such that for each
connected open neighborhood $V_{F}\subset V_{0; F}$, we get a unique, continuously bounded differentiable
map $T\colon V_{F}\to \R$ with
\begin{equation} \label{Tabdiff2}
T(F)=x_0, \quad f'(T(f))=0,\;\;\; f\in V_{F},\quad dT(f)h=-h'(T(f))/f''(T(f))
\end{equation}
But by assumption~\eqref{nahbei} from some $n$ on, $F_n$ and $G_n$ will lie in $V_{0; F}$, and setting $x_n=T(F_n)$,
by \eqref{Tabdiff2} $F_n'(x_n)=0$, and
$$|x_n-x_0| =|T(F_n)-T(F)|\leq |F_n'(x_0)|/F''(x_0)=\LO(n^{-\beta'})$$
which is (b); again by \eqref{nahbei},
\begin{eqnarray*}
|F_n''(x_n)-F''(x_0)|&\leq&  |F_n''(x_n)-F''(x_n)| +|F''(x_n)-F''(x_0)| \leq 
\sup_x |F_n''(x)-F''(x)| +\Lo(n^0) =\Lo(n^0)
\end{eqnarray*}
In particular, eventually in $n$, $F_n''(x_n)>0$ and hence $x_n$ is a minimum of $F$, so (a) is shown.
By \eqref{nahbei},
$
\sup_x |F-G_n| + |F'-G'_n| + |F''-G''_n|=\LO(n^{-\beta'})
$,
so (c) follows just as (a). For (d) we note
$$|x_n-y_n| =|T(F_n)-T(G_n)|\leq |G_n'(x_n)|/F_n''(x_n)\stackrel{(a)}=|G_n'(x_n)|/(f_2+\Lo(n^0))=\LO(n^{-\beta})$$
To show (e),  we introduce $d_n:=y_n-x_n$ and write
\begin{equation}
0\leq G_n(x_n)-G_n(y_n)=G_n'(y_n)d_n+G_n''(y_n)d_n^2/2+\Lo(d_n^2)=
(f_2+\Lo(n^0))d_n^2/2+\Lo(d_n^2)=\LO(n^{-2\beta})
\end{equation}
\hfill\qed

\subsection[Proof of Proposition~xxx]{Proof of Proposition~\ref{sattel}}\label{psattel}
We show that under the assumptions of this proposition
$x_j$  indeed defines a ``uniformly bad contamination'' in the sense that
for the fixed contamination $Q_n(x_j)$
\begin{equation}
 {\rm asMSE}_0(S^{\SSs (b_0)}_n,Q_n(x_0)) =\min_{b>0} {\rm asMSE}_0(S^{\SSs (b)}_n,Q_n(x_0))
\end{equation}
resp.\
$
 {\rm asMSE}_1(S^{\SSs (c_1)}_n,Q_n(x_1)) =\min_{c>0} {\rm asMSE}_1(S^{\SSs (c)}_n,Q_n(x_1))
$ 
In case $j=0$, as in the setup of \citet[chap.~5]{Ri:94}, we obtain 
\begin{equation}
{\rm asMSE}_0(S^{\SSs (b)}_n,Q_n(x_0))=\tr \Cov_{\rm \SSs id}(\eta_b)+r^2|\eta_b(x_0)|^2,\qquad
{\rm asMSE}_0(\hat S_n,Q_n(x_0))=\tr {\cal I} +r^2|\hat \eta(x_0)|^2
\end{equation}
Now for given $x_0$, either $|\eta^{(b)}(x_0)|<b$ or $|\eta^{(b)}(x_0)|=b$. In the first case, \eqref{minetahat} applies and hence
\begin{equation}
  {\rm asMSE}_0(S^{\SSs (b_0)}_n,Q_n(x_0))\geq {\rm asMSE}_0(\hat S_n,Q_n(x_0))
\end{equation}
In the latter, $Q_n(x_0)$ already achieves maximal \asy risk for $S^{\SSs (b)}_n$ on ${\cal Q}_n$, and hence
by minimaxity of $S^{\SSs (b_0)}_n$
\begin{equation}
{\rm asMSE}_0(S^{\SSs (b)}_n,Q_n(x_0))\geq {\rm asMSE}_0(S^{\SSs (b_0)}_n,Q_n(x_0))
\end{equation}
For the case $j=1$ one argues in an analogue way.
\hfill\qed

\subsection[Proof for (xx) in the Gaussian location scale model]{Proof for \eqref{minetahat} in the Gaussian location scale model}\label{minetp}
    We abbreviate the location and scale parts by indices $l$ and $s$ respectively.
    By equivariance we may limit ourselves to the case $\theta=(0,1)^{\tau}$. Due to symmetry, $A=A(b)$ from \eqref{allgemoptIC} is diagonal for all $b$
    with elements $A_{l}$ and $A_s$ and we may write
    \begin{equation}
\eta_b=Y\min\{1,b/|Y|\},\qquad Y^\tau =\big(A_lx,A_s(x^2-1-z_s)\big)
    \end{equation}
    The centering $z_s(b)$ after the clipping is necessary, as the scale part is not skew symmetric;
    in the pure scale case (with known $\theta_l$), the corresponding centering $z'_s=z'_s(b)$ is antitone in $b$, because $\Lambda_s$ is
    monotone in $x^2$:
    It decreases from $0$ to $[\Phi^{-1}(3/4)]^2-1\doteq -0.545=:\check z$. In the combined case,
    we never reach this extremal case due to the additional location part---compare \citet[Remark~8.2.1(a)]{Koh:04di}
    where $\bar z_s=\bar a_{\rm\SSs sc}/\bar \alpha -1\doteq -0.530$;  in any case, $z_s>-1$ always.
    Hence in particular, for $x_0=1.844$ and $b$ such that $|\eta^{(b)}(x_0)|\leq b$ it holds that
\begin{equation}
|\eta^{(b)}_{s}(x_0)|=A_{s}(b)|x_0^2-1-z_{s}(b)|>A_{s}(b)|x_0^2-1|>
{\cal I}_{s}^{-1}|x_0^2-1|
=|\hat \eta_{s}(x_0)|
\end{equation}
and thus in particular,
\begin{eqnarray}
|\eta^{(b)}(x_0)|^2&=&|\eta^{(b)}_{s}(x_0)|^2+
|\eta^{(b)}_{l}(x_0)|^2)=|\eta^{(b)}_{s}(x_0)|^2+
A_{0;l}(b) x_0^2>
\hat \eta_{s}(x_0)^2+
{\cal I}_{l}^{-2}x_0^2
=|\hat \eta_{}(x_0)|^2
\end{eqnarray}
\hfill\qed

\section*{Acknowledgement}
\medskip

\medskip
\hrulefill\hspace*{6cm}\\[2ex]
Web-page to this article:\\
\href{http://www.mathematik.uni-kl.de/~ruckdesc/}%
{{\footnotesize \url{http://www.mathematik.uni-kl.de/~ruckdesc/}}}

%
%
%
%
%
\end{document}